\documentclass[12pt]{article}
\usepackage[T1]{fontenc}

\usepackage{amsfonts}
\usepackage{amssymb}
\usepackage{latexsym}
\usepackage{color}
\usepackage{graphicx}
\usepackage{algorithmic}

\newtheorem{theorem}{Theorem}
\newtheorem{corollary}{Corollary}
\newtheorem{lemma}{Lemma}

\newtheorem{example}{Example}

\newtheorem{definition}{Definition}
\newtheorem{remark}{Remark}

\newcommand{\bfe}{\mbox{$\mbox{\boldmath $e$}$}} 
\newcommand{\bff}{\mbox{$\mbox{\boldmath $f$}$}} 
\newcommand{\bfp}{\mbox{$\mbox{\boldmath $p$}$}} 
\newcommand{\bfq}{\mbox{$\mbox{\boldmath $q$}$}} 
\newcommand{\bfs}{\mbox{$\mbox{\boldmath $s$}$}} 
\newcommand{\bfu}{\mbox{$\mbox{\boldmath $u$}$}} 
\newcommand{\bfx}{\mbox{$\mbox{\boldmath $x$}$}} 
\newcommand{\bfy}{\mbox{$\mbox{\boldmath $y$}$}} 

\newcommand{\bxi}{\mbox{$\mbox{\boldmath $\xi$}$}} 

\newcommand{\bfeta}{\mbox{$\mbox{\boldmath $\eta$}$}} 
\newcommand{\blambda}{\mbox{$\mbox{\boldmath $\lambda$}$}}

\newcommand{\bfzero}{\mbox{$\mbox{\boldmath $0$}$}} 

\date{}

\title{Inversion of operator pencils on Banach space using Jordan chains when the generalized resolvent has an isolated essential singularity}

\author{Amie Albrecht, Phil Howlett, Geetika Verma}

\begin{document}

\maketitle

\begin{abstract}
We assume that the generalized resolvent for a bounded linear operator pencil mapping one Banach space onto another has an isolated essential singularity at the origin and is analytic on some annular region of the complex plane centred at the origin. In such cases the resolvent operator can be represented on the annulus by a convergent Laurent series and the spectral set has two components\textemdash a bounded component inside the inner boundary of the annulus and an unbounded component outside the outer boundary.  In this paper we prove that the complementary spectral separation projections on the domain space are uniquely determined by the respective generating subspaces for the associated infinite-length generalized Jordan chains of vectors and that the domain space is the direct sum of these two subspaces.  We show that the images of the generating subspaces under the mapping defined by the pencil provide a corresponding direct sum decomposition for the range space and that this is simply the decomposition defined by the complementary spectral separation projections on the range space.  If the domain space has a Schauder basis we show that the separated systems of fundamental equations are reduced to two semi-infinite systems of matrix equations which can be solved recursively to obtain a basic solution and thereby determine the Laurent series coefficients for the resolvent operator on the given annular region.  
\end{abstract}



\section{Introduction}
\label{intro}

Let $X,Y$ be complex Banach spaces and let $A_0, A_1 \in {\mathcal B}(X,Y)$, where $A_0$ is singular.  Define a linear operator pencil $A: {\mathbb C} \rightarrow {\mathcal B}(X,Y)$ by the formula $A(z) = A_0 + A_1z$.  It was shown by Albrecht et al.~\cite{alb2} that the resolvent $R:\mathcal U_{s,r} \rightarrow {\mathcal B}(Y,X)$ for $A$ on the annulus ${\mathcal U}_{s,r} = \{z \in {\mathbb C} \mid s < |z| < r\}$ is defined by a Laurent series $R(z) = \sum_{j \in {\mathbb Z}} R_jz^j$  if and only if the coefficients $R_j \in {\mathcal B}(Y,X)$ for each $j \in {\mathbb Z}$ satisfy a system of left and right fundamental equations and are suitably bounded.  If a solution exists then it is uniquely defined by a basic solution $\{R_{-1},R_0\}$ and there is a closed form $R(z) = (I z + R_{-1}A_0)^{-1}R_{-1} + (I + R_0A_1z)^{-1}R_0$ for the resolvent.  Furthermore the functions ${\mathcal R}_{\lambda} = \lambda^{-1}R(-\lambda^{-1})A_0 \in {\mathcal B}(X)$ and ${\mathcal S}_{\lambda} = \lambda^{-1}A_0R(-\lambda^{-1}) \in {\mathcal B}(Y)$ each satisfy a classical resolvent equation.  The fundamental equations were also used in \cite{alb2} to show that the operators $P = R_{-1}A_1 \in {\mathcal B}(X)$ and $Q = A_1R_{-1}\in {\mathcal B}(Y)$ define corresponding key projections that separate the bounded and unbounded components $\sigma_1 = \{ z \in {\mathbb C} \mid |z| \leq s \}$ and $\sigma_2 = \{ z \in {\mathbb C} \mid |z| \geq r \}$ of the spectral set $\sigma$ and that the key projections can be used to establish useful properties of the coefficients $\{R_j\}_{j \in {\mathbb Z}}$.  Albrecht et al.~\cite{alb2} extended these results to polynomial pencils and investigated the global structure of $R(z)$ in the case where $R(z)$ has only a finite number of isolated singularities.

We wish to find a solution procedure for the fundamental equations.  We will show that the complementary projections on the domain space which separate the two components of the spectral set on the given annulus\textemdash the inner bounded component and the outer unbounded component\textemdash are uniquely determined by the respective generating subspaces for the associated infinite-length generalized Jordan chains of vectors.  We will use these chains to find the key projection operators and the desired direct sum decomposition for the domain and range spaces.  We can then find a basic solution to the fundamental equations and hence construct the Laurent series representation for the generalized resolvent.  We motivate our discussion by showing that the necessary and sufficient conditions obtained by Bart and Lay \cite{bar1} when the resolvent has a finite-order pole\textemdash that the ascent and descent of the unperturbed operator with respect to the perturbation operator are equal to the order of the pole\textemdash have no direct analogy at an isolated essential singularity.  Finally we will use a typical linear operator pencil on a standard separable Banach space to demonstrate the new results and show how they can be applied to find the resolvent operator.

\section{Notation}
\label{not}

For the most part we use standard notation.  However if $E$ is a subset of a Banach space $X$ we follow Yosida \cite{yos1} and write $E^a$ for the closure of the set $E$.  Let $X,Y$ be Banach spaces.  If $T \subseteq Y$ is a subspace and $A \in {\mathcal B}(X,Y)$ is a bounded linear operator then it is standard notation to write $A^{-1}(T) = \{\bfx \in X \mid A \bfx \in T\}$ whether or not $A$ is invertible.  We will use an extended form of this notation.  Let $A, B \in {\mathcal B}(X,Y)$.  If $S_0 \subseteq X$ is a subspace we will write $S_n = (A^{-1}B)^n(S_0) = A^{-1}B(S_{n-1}) = \{ \bfx \in X \mid A \bfx \in B(S_{n-1}) \}$ for each $n \in {\mathbb N}$ whether or not $A$ is invertible.  For $s \in {\mathbb R}$ and $r \in {\mathbb R} \cup \{\infty\}$ with $0 \leq s < r \leq \infty$ we will make frequent use of the notation ${\mathcal U}_{s,r} = \{ z \in {\mathbb C} \mid s < |z| < r\}$ to denote an open annular subset of the complex plane ${\mathbb C}$ centred at $z = 0$ with inner radius $s$ and outer radius $r$.

\section{The main result}
\label{mt}

In order to state our main result we need to preview some later definitions of the infinite-length generalized Jordan chains.  The sequence $\{ \bfx_{-n} \}_{n \in {\mathbb N}} \subseteq X$ is an infinite-length {\em singular} Jordan chain generated by $\bfx_{-1}$ for the pencil $A(z)$ on the annular region ${\mathcal U}_{s,\, \infty}$ if $A_0 \bfx_{-n} + A_1 \bfx_{-n-1} = \bfzero$ for all $n \in {\mathbb N}$ and $\lim_{n \rightarrow \infty} \| \bfx_{-n} \|^{1/n} = s$.  The sequence $\{ \bfx_{n} \}_{n \in {\mathbb N}} \subseteq X$ is an infinite-length {\em regular} Jordan chain generated by $\bfx_1$ for the pencil $A(z)$ on the annular region ${\mathcal U}_{0,\,r}$ if $A_0 \bfx_{n} + A_1 \bfx_{n-1} = \bfzero$ for all $n \in {\mathbb N}$ and $\lim_{n \rightarrow \infty} \| \bfx_{n} \|^{1/n} = 1/r$.

\begin{theorem}
\label{th1}
{\rm Let $0 \leq s < r \leq \infty$ and suppose that $R(z)$ is analytic for $z \in {\mathcal U}_{s,r}$.  Let $P = R_{-1}A_1 \in {\mathcal B}(X)$ and $P^c = I - P = R_0A_0 \in {\mathcal B}(X)$ be the complementary key projections on the domain space that separate the bounded and unbounded parts of the spectral set for $A(z)$ relative to the annular region ${\mathcal U}_{s,r}$ and let $Q = A_1R_{-1} \in {\mathcal B}(Y)$ and $Q^c = I - Q = A_0R_0 \in {\mathcal B}(Y)$ be the corresponding complementary key projections on the range space.  Now let $X_s \subseteq X$ and $X_r \subseteq X$ be the respective generating subspaces for the infinite-length singular and regular generalized Jordan chains for $A(z)$ on ${\mathcal U}_{s,\, \infty}$ and ${\mathcal U}_{0,\, r}$.  Then we have (i) $X_s = P(X)$ and $X_r = P^c(X)$ with $X = X_s \oplus X_r \cong X_s \times X_r$ and (ii) $Y_s = Q(Y) = A_1(X_s)$ and $Y_r = Q^c(Y) = A_0(X_r)$ with $Y = Y_s \oplus Y_r \cong Y_s \times Y_r$.} $ \hfill \Box$
\end{theorem}

This theorem will be justified in Section \ref{ilgjc} by establishing certain definitive properties for the infinite-length generalized Jordan chains and the associated generating subspaces.

\section{Previous work}
\label{pw}

We restrict our detailed discussion to the key papers by Stummel \cite{stu1}, Bart and Lay \cite{bar1} and Albrecht et al.~\cite{alb2}.  We will not describe other recent contributions by Albrecht et al.~\cite{alb1, alb3, alb4}, Avrachenkov et al.~\cite{avr2,avr3}, Franchi \cite{fra1} and Howlett et al.~\cite{how2,how3} nor the earlier work by Howlett \cite{how1}, Langenhop \cite{lan1,lan2}, Rose \cite{ros1}, Rothblum \cite{rot1}, Schweitzer and Stewart \cite{sch1} and Wilkening \cite{wil1} but rather refer readers to the more expansive reviews in our two recent papers \cite{alb3,alb4}.  A complete discussion about spectral methods for operator pencils can be found in the book by Gohberg et al.~\cite[pp 49\textendash 54]{goh1} and we refer to the classical texts by Kato \cite[pp 178\textendash 179]{kat1} and Yosida \cite[pp 228\textendash 231]{yos1} for further background information.

The first major contribution to a general theory of inversion for operator pencils on Banach space was a comprehensive paper by F. Stummel \cite{stu1}.  The paper used complex function theory and contour integrals to define key projections that separate the singular and regular components of the generalized resolvent operator.  Although the work by Stummel provides a complete structural answer to the problem the contour integrals cannot be used to calculate the key projections or to calculate the Laurent series coefficients if the resolvent is not known {\em a priori}.  This raises a legitimate question about systematic calculation of the generalized resolvent for operator pencils on Banach space.

An elegant paper by Bart and Lay \cite{bar1} used generalized Jordan chains of subspaces to answer the question of systematic calculation in the case where the generalized resolvent has a finite-order pole.  This work extends the method developed for matrix pencils by Langenhop \cite{lan1,lan2} to pencils of linear operators $A(z) = A_0 + A_1z \in {\mathcal B}(X,Y)$ where $X, Y$ are complex Banach spaces.  For each such pencil Bart and Lay define corresponding decreasing sequences of subspaces $\{S_m\}_{m \in {\mathbb N}-1} \subseteq X$ and $\{T_m\}_{m\in {\mathbb N}-1} \subseteq Y$ where $S_m = (A_1^{-1}A_0)^m(X) \subseteq S_0 = X$ and $T_m = A_0(A_1^{-1}A_0)^{m-1}(Y) \subseteq T_0 = Y$ and corresponding increasing sequences of subspaces $\{U_m\}_{m \in {\mathbb N}-1} \subseteq X$ and $\{V_m\}_{m\in {\mathbb N}-1} \subseteq Y$ where $\{ \bfzero\} = U_0 \subseteq U_m = A_0^{-1}(A_1A_0^{-1})^{m-1}(\{ \bfzero\}) \subseteq X$ and $\{ \bfzero \} = V_0 \subseteq V_m = (A_1A_0^{-1})^m(\{ \bfzero\}) \subseteq Y$ for each $m \in {\mathbb N}$.  The ascent $\alpha(A_0:A_1)$ of $A_0$ relative to $A_1$ is the smallest extended natural number $m \in {\mathbb N} \cup \{\infty\}$ with $A_0^{-1}(\{\bfzero\}) \cap S_m = \{ \bfzero\}$ and the descent $\delta(A_0:A_1)$ of $A_0$ relative to $A_1$ is the smallest extended natural number $m \in {\mathbb N} \cup \{\infty\}$ with $A_0(X) \oplus V_m = Y$.  The subspaces defined above were used to show that if the ascent and descent are both finite then they are equal and that $\alpha(A_0:A_1) = \delta(A_0:A_1) = m < \infty$ if and only if $R(z) = A(z)^{-1}$ is analytic on a region ${\mathcal U}_{0,r}$ with a pole of order $m$ at $z=0$.  Bart and Lay did not consider isolated essential singularities.

An alternative approach to the same problem was provided by Albrecht et al.~\cite{alb2} who used a system of fundamental equations to define the generalized resolvent for the pencil $A(z) = A_0 + A_1z \in {\mathcal B}(X,Y)$ where $X, Y$ are complex Banach spaces.  They showed that the resolvent $R(z) = A(z)^{-1}$ is analytic on an annular region ${\mathcal U}_{s,r}$ where $0 \leq s < r \leq \infty$ and can be represented by a Laurent series $R(z) = \sum_{j \in {\mathbb Z}} R_j z^j$ where $R_j \in {\mathcal B}(Y,X)$ for all $j \in {\mathbb Z}$ if and only if
\begin{equation}
\label{rcb}
\lim_{k \rightarrow \infty} \|R_{-k}\|^{1/k} \leq s\ \mbox{and}\ \lim_{\ell \rightarrow \infty} \|R_{\ell}\|^{1/\ell} \leq 1/r
\end{equation}
and the coefficients $\{R_j\}_{j \in {\mathbb Z}}$ satisfy the left fundamental equations
\begin{equation}
\label{lfe}
R_{j-1}A_1 + R_jA_0 = \left \{ \begin{array}{ll}
I & \mbox{if}\ j=0 \\
0 & \mbox{if}\ j \neq 0 \end{array} \right.
\end{equation}
and the right fundamental equations
\begin{equation}
\label{rfe}
 A_1R_{j-1} + A_0R_j = \left \{ \begin{array}{ll}
I & \mbox{if}\ j=0 \\
0 & \mbox{if}\ j \neq 0. \end{array} \right.
\end{equation}

More precisely Albrecht et al.~\cite{alb2} proved the following result.

\begin{theorem}
\label{th2}
{\rm The coefficients $\{R_j\}_{j \in {\mathbb Z}} \in {\mathcal B}(Y,X)$ satisfy $(\ref{rcb})$, $(\ref{lfe})$ and $(\ref{rfe})$ if and only if the following are all satisfied: $(i)$ $P = R_{-1}A_1 \in {\mathcal B}(X)$ and $P^c = I-P = R_0A_0 \in {\mathcal B}(X)$ are complementary projections on $X$; and $Q = A_1R_{-1} \in {\mathcal B}(Y)$ and $Q^c = I - Q = A_0R_0 \in {\mathcal B}(Y)$ are corresponding complementary projections on $Y$;  $(ii)$ $A_i = QA_iP + Q^cA_iP^c$ for $i=0,1$;  $(iii)$ $R_j = PR_jQ$ for $j \leq -1$ and $R_j = P^cR_jQ^c$ for $j \geq 0$; $(iv)$ $R_{-k} = (-1)^{k-1}(R_{-1}A_0)^{k-1}R_{-1}$ for $k \in {\mathbb N}$ and $R_{\ell} = (-1)^{\ell} (R_0A_1)^{\ell}R_0$ for $\ell \in {\mathbb N} - 1$; and $(v)$ $\lim_{k \rightarrow \infty} \|(R_{-1}A_0)^k\|^{1/k} \leq s$  and $\lim_{\ell \rightarrow \infty} \|(R_0A_1)^{\ell} \|^{1/{\ell}} \leq 1/r$.} $ \hfill \Box$
\end{theorem}

The above result suggests that the coefficients $\{R_j\}_{j \in {\mathbb Z}}$ could be found directly by solving a system of fundamental equations without prior knowledge of a general expression for $R(z)$.  In this regard it has been shown \cite[Corollary 2]{alb2} that the resolvent $R:{\mathcal U}_{s,r} \rightarrow {\mathcal B}(Y,X)$ is analytic if and only if there exists a basic solution $\{R_{-1}, R_0\} \in {\mathcal B}(Y,X)$ such that $(i)$ $R_{-1}A_1 + R_0A_0 = I$ and $A_1R_{-1} + A_0R_0 = I$;  $(ii)$ $R_{-1}A_iR_0 = 0$ and $R_0A_iR_{-1} = 0$ for each $i = 0,1$; and $(iii)$ $\lim_{k \rightarrow \infty}  \|(R_{-1}A_0)^k\|^{1/k} \leq s$ and $\lim_{\ell \rightarrow \infty} \|(R_0A_1)^{\ell} \|^{1/{\ell}} \leq 1/r$ and that in this case $R_{-k} = (-1)^{k-1}(R_{-1}A_0)^{k-1}R_{-1}$ for $k \in {\mathbb N}$ and $R_{\ell} = (-1)^{\ell} (R_0A_1)^{\ell}R_0$ for $\ell \in {\mathbb N}-1$.  This latter result is an analogue of the original matrix result established by Langenhop in \cite{lan2}.

\section{Two related examples where the resolvent has an isolated essential singularity}
\label{exries}

Let $X = Y = \ell^2$ be the Hilbert space of square-summable sequences of complex numbers and consider a linear operator pencil $A(z) = A_0 + A_1z$ where $A_0, A_1 \in {\mathcal B}(\ell^2)$ are bounded linear operators.  We will use two related examples to show that the necessary and sufficient conditions established by Bart and Lay \cite{bar1} for finite-order poles of the resolvent do not extend directly to isolated essential singularities.

\begin{example}
\label{ex1}

{\rm Let $A_0 = B$ where
$$
B = \left[ \begin{array}{cccccc}
0 & \beta_1 & 0 & 0 & 0 & \cdots \\
0 & 0 & \beta_2 & 0 & 0 & \cdots \\
0 & 0 & 0 & \beta_3 & 0 & \cdots \\
0 & 0 & 0 & 0 & \beta_4 & \cdots \\
\vdots & \vdots & \vdots & \vdots & \vdots & \ddots \end{array} \right] \in {\mathcal B}(\ell^2)
$$
and $\{ \beta_n \}_{n \in {\mathbb N}} \in {\mathbb C}$ is a known sequence of complex numbers with $|\beta_n|^{1/n} \rightarrow 0$ as $n \rightarrow \infty$ and let $A_1 = I$ where $I \in {\mathcal B}(\ell^2)$ denotes the identity operator.  The resolvent $R(z) = A(z)^{-1} \in {\mathcal B}(\ell^2)$ is well defined by the Neumann expansion $R(z) = I/z - B/z^2 + B^2/z^3 - B^3/z^4 + \cdots$ for all $z \in {\mathcal U}_{\, 0,\infty} = \{z \in {\mathbb C} \mid 0 < |z| < \infty \}$.  If we write $B^m = [b_{i,j}^{(m)}]$ in matrix form for each $m \in {\mathbb N}$ then
$$
b_{i,j}^{(m)} = \left\{ \begin{array}{ll}
\prod_{k=i}^{i+m-1} \beta_k & \mbox{for}\ i \in {\mathbb N}\ \mbox{and}\ j = i+m \\
& \\
0 & \mbox{otherwise}. \end{array} \right.
$$
In this case the analogous Jordan chains of subspaces to those defined by Bart and Lay \cite{bar1} are given by $S_0 = T_0 = \ell^2$ with $S_m = T_m = [B^m(\ell^2)]^a$ for each $m \in {\mathbb N}$ where $[B^m(\ell^2)]^a$ denotes the closure of $B^m(\ell^2)$ and $U_0 = V_0 = \{ \bfzero\}$ with $U_m = V_m = B^{-m}(\{ \bfzero\}) = \{ \bfx \in \ell^2 \mid B^m \bfx = \bfzero \}$ for each $m \in {\mathbb N}$.  

Let $\{\bfe_j\}_{j \in {\mathbb N}}$ denote the standard orthonormal basis in $\ell^2$.  If $\bfy = y_1 \bfe_1 + y_2 \bfe_2 + \cdots \in {\ell}^2$ then for each $\epsilon > 0$ we can find $p = p(\epsilon) \in {\mathbb N}$ such that $\bfy^{(p)} = y_1\bfe_1 + \cdots + y_p \bfe_p$ satisfies $\| \bfy - \bfy_p \| < \epsilon$.   If we choose $\bfx^{(p)} = [x_i]$ such that 
$$
x_{m+i} = \left\{ \begin{array}{ll}
y_{i}/b_{i,m+i}^{(m)} & \mbox{for}\ i = 1, \ldots, p \\
& \\
0 & \mbox{otherwise} \end{array} \right.
$$
then $B^m \bfx^{(p)} = \bfy^{(p)}$.  Therefore $S_m = T_m = [B^m(\ell^2)]^a = \ell^2$ for each $m \in {\mathbb N}$.  The equation $B^m \bfx = \bfzero$ has solutions $\bfx = [x_j] \in {\ell^2}$ given by 
$$
x_j = \left\{ \begin{array}{ll}
t_j & \mbox{for}\ j=1,\ldots,m \\
& \\
0 & \mbox{otherwise} \end{array} \right.
$$
where $t_1,\ldots,t_m \in {\mathbb C}$ are arbitrary parameters.  It follows that $U_m = V_m = \{ \bfx \in {\ell}^2 \mid \bfx  = t_1 \bfe_1 + \cdots + t_m \bfe_m\}$ for each $m \in {\mathbb N}$.  We follow Bart and Lay \cite{bar1} and define the ascent of $B$ relative to $I$ by setting
$$
\alpha(B:I) = \inf \{m \in {\mathbb N} \cup \infty \mid U_1 \cap S_m = \{ \bfzero\} \}
$$
and the descent of $B$ relative to $I$ by setting
$$
\delta(B:I) = \inf \{m \in {\mathbb N} \cup \infty \mid T_1 + V_m = \ell^2 \}.
$$
Since $U_1 \cap S_m = \{ t_1 \bfe_1 \} \cap \ell^2 = \{ t_1 \bfe_1 \} \neq \{ \bfzero \}$ for all $m \in {\mathbb N} - 1$ we have $\alpha(B:I) = \infty$ and since $T_1 + V_0 = \ell^2 + \{ \bfzero \} = \ell^2$ we have $\delta(B:I) = 0$.   Thus the resolvent has an essential singularity at $z =0$ with $0 = \delta(B:I) \neq \alpha(B:I) = \infty$.}  $\hfill \Box$
\end{example}

\begin{example}
\label{ex2}

{\rm Let $A_0 = C$ where
$$
C = \left[ \begin{array}{ccccc}
0 & 0 & 0 & 0 & \cdots \\
\gamma_1 & 0 & 0 & 0 & \cdots \\
0 & \gamma_2 & 0 & 0 & \cdots \\
0 & 0 & \gamma_3 & 0 & \cdots \\
0 & 0 & 0 & \gamma_4 & \cdots \\
\vdots & \vdots & \vdots & \vdots & \ddots \end{array} \right] 
$$
and $\{ \gamma_n \}_{n \in {\mathbb N}} \in {\mathbb C}$ is a known sequence of complex numbers with $|\gamma_n|^{1/n} \rightarrow 0$ as $n \rightarrow \infty$ and let $A_1 = I$ where $I \in {\mathcal B}(\ell^2)$ denotes the identity operator.  The resolvent $R(z) = A(z)^{-1} \in {\mathcal B}(\ell^2)$ is well defined by the Neumann expansion $R(z) = I/z - C/z^2 + C^2/z^3 - C^3/z^4 + \cdots$ for all $z \in {\mathcal U}_{\, 0,\infty} = \{z \in {\mathbb C} \mid 0 < |z| < \infty \}$.  Similar calculations to those used in Example \ref{ex1} show that  the resolvent has an essential singularity at $z =0$ with $0 = \alpha(C:I) \neq \delta(C:I) = \infty$.}  $\hfill \Box$
\end{example}

\section{The infinite-length generalized Jordan chains}
\label{ilgjc}

Let $X$ and $Y$ be complex Banach spaces.  We assume that the resolvent $R(z) \in {\mathcal B}(Y,X)$ of the linear pencil $A(z) = A_0 + A_1z \in {\mathcal B}(X,Y)$ is analytic on some annular region ${\mathcal U}_{s,r} = \{ z \in {\mathbb C} \mid s < |z| < r \}$ where $0 \leq s < r \leq \infty$.   Thus there exist operator coefficients $\{R_j\}_{j \in {\mathbb Z}} \in {\mathcal B}(Y,X)$ satisfying the properties described in Theorem \ref{th2} with $R(z) = \sum_{j \in {\mathbb Z}} R_j z^j$ for all $z \in {\mathcal U}_{s,r}$.  Our investigation will focus on construction of the key projection operators $P = R_{-1}A_1 \in {\mathcal B}(X)$ and $Q = A_1R_{-1} \in {\mathcal B}(Y)$ that separate the bounded ($|z| \leq s$) and unbounded ($|z| \geq r$) components of the spectral set for the pencil $A(z)$.

\begin{remark}
\label{rem1}
{\rm  Although our assumption that $R(z) \in {\mathcal B}(Y,X)$ is analytic for $z \in {\mathcal U}_{\, s,r}$ guarantees existence of the coefficients $\{ R_j\}_{j \in {\mathbb Z}} \in {\mathcal B}(Y,X)$ and although we use these coefficients in our subsequent arguments to establish the key structural properties we emphasize that the coefficients are not known {\em a priori}.  Our aim is to find a systematic methodology to determine the coefficients given that they exist.} $\hfill \Box$
\end{remark}

Theorem \ref{th1} is a consequence of the definitive properties of the infinite-length generalized Jordan chains, which we establish in this section.  These results are new.  The properties of the finite-length Jordan chains are reformulated from the theory developed by Bart and Lay \cite{bar1}, however the definition of singular and regular chains is new.

\subsection{The finite-length singular Jordan chains}
\label{flsjc}

To begin we follow the methodology of Bart and Lay \cite{bar1} and define subspaces $U_{\mbox{\scriptsize sg},\,0} = \{ \bfzero \} \subseteq X$ and
$$
U_{\mbox{\scriptsize sg},\, k} = A_0^{-1}(A_1A_0^{-1})^{k-1}( \{ \bfzero \}) \subseteq X
$$
for each $k \in {\mathbb N}$.  We also define corresponding subspaces $V_{\mbox{\scriptsize sg},\,0} = \{ \bfzero \} \subseteq Y$ and
$$
V_{\mbox{\scriptsize sg},\,k} = (A_1 A_0^{-1})^k( \{ \bfzero \}) \subseteq Y
$$
for each $k \in {\mathbb N}$.  We will use these subspaces to define the finite-length singular Jordan chains of vectors for the pencil $A(z) = A_0 + A_1z$.  We observe that $\bfx_{-1} \in U_{\mbox{\scriptsize sg},\,k} \Leftrightarrow A_0 \bfx_{-1} \in V_{\mbox{\scriptsize sg},\,k-1} = A_1(U_{\mbox{\scriptsize sg},\,k-1})$.  Thus there exists $\bfx_{-2} \in U_{\mbox{\scriptsize sg},\,k-1}$ such that $A_0 \bfx_{-1} = A_1 (- \bfx_{-2}) \Leftrightarrow A_0 \bfx_{-1} + A_1 \bfx_{-2} = \bfzero$.  By repeating this argument it follows that we can find a Jordan chain $\{ \bfx_{-n} \}_{n=1}^{k}$ with $\bfx_{-n} \in U_{\mbox{\scriptsize sg},\,k-n+1}$ for each $n=1,\ldots,k$ such that
\begin{equation}
\label{ens}
A_0 \bfx_{-n} + A_1 \bfx_{-n-1} = \bfzero
\end{equation}
for each $n = 1,\ldots,k-1$.  Since $\bfx_{-k} \in U_{\mbox{\scriptsize sg},\,1} = A_0^{-1}( \{ \bfzero\})$ it follows that
\begin{equation}
\label{eks}
A_0 \bfx_{-k} = \bfzero.
\end{equation}
We say that $\{ \bfx_{-n} \}_{n=1}^{k}$ is a {\em singular} Jordan chain of length $k$ for $A(z)$ that is generated by the point $\bfx_{-1} \in U_{\mbox{\scriptsize sg},\,k}$.  The subspace $U_{\mbox{\scriptsize sg},\,k}$ is the generating subspace for all singular Jordan chains of length $k$ for the pencil $A(z)$.  We say these chains are singular because the generating set is closely related to the coefficients of the negative powers in the Laurent series expansion for $R(z) = A(z)^{-1}$ on the region $z \in {\mathcal U}_{\, s,r}$.  Note that $A_0^{-1}( \{ \bfzero\}) = U_{\mbox{\scriptsize sg},\,1} \subseteq U_{\mbox{\scriptsize sg},\,2} \subseteq \cdots \subseteq U_{\mbox{\scriptsize sg},\,k} \subseteq X$.   We have established the following result.

\begin{lemma}
\label{lem1}
{\rm  For each $k \in {\mathbb N}$ the subspace $U_{\mbox{\scriptsize sg},\,k} = A_0^{-1}(A_1A_0^{-1})^{k-1}( \{ \bfzero \}) \subseteq X$ is the generating subspace for all singular Jordan chains of length $k$ for the pencil $A(z)$.  Thus we have 
$$
U_{\mbox{\scriptsize sg},\,k} = \{ \bfx_{-1} \in X \mid\ \mbox{there exists}\ \{ \bfx_{-n} \}_{n=1}^k\ \mbox{with}\ \bfx_{-n} \in U_{\mbox{\scriptsize sg},\,k-n+1}\ \mbox{such that (\ref{ens}) and (\ref{eks}) are satisfied} \} 
$$
for each $k \in {\mathbb N}$. $\hfill \Box$}
\end{lemma}

Consider a singular Jordan chain $\{\bfx_{-n}\}_{n=1}^k$ of length $k$ generated by a point $\bfx_{-1} \in U_{\mbox{\scriptsize sg},\,k}$.  Now
$$
\bfx_{-k} \in U_{\mbox{\scriptsize sg},\,1} \Rightarrow A_0 \bfx_{-k} = \bfzero \Rightarrow R_0A_0 \bfx_{-k} = \bfzero \Leftrightarrow P^c\bfx_{-k} = \bfzero \Leftrightarrow P \bfx_{-k} = \bfx_{-k}.
$$
Suppose $P\bfx_{-k+n-1} = \bfx_{-k+n-1} \Leftrightarrow R_{-1}A_1 \bfx_{-k+n-1} = \bfx_{-k+n-1}$ for some $n = m \in {\mathbb N}$ with $1 \leq m < k$.  The sequence $\{\bfx_{-n}\}_{n=1}^k$ is a singular Jordan chain and so $A_0 \bfx_{-k+m} + A_1 \bfx_{-k+m-1} = \bfzero$.  Therefore $R_{-1}A_0 \bfx_{-k+m} + R_{-1}A_1 \bfx_{-k+m-1} = \bfzero$ and our supposition that $R_{-1}A_1 \bfx_{-k+m-1} = \bfx_{-k+m-1}$ now means that
$$
\bfx_{-k+m-1} = (-1) R_{-1}A_0 \bfx_{-k+m}.
$$
The relationship $A_0 \bfx_{-k+m} + A_1 \bfx_{-k+m-1} = \bfzero$ also implies $R_0A_0 \bfx_{-k+m} + R_0 A_1 \bfx_{-k+m-1} = \bfzero$ which in turn gives
$$
R_0 A_0 \bfx_{-k+m} = (-1) R_0A_1 \bfx_{-k+m-1} = (-1)^2 R_0A_1R_{-1}A_1 \bfx_{-k+m-1}.
$$
By applying the identity $R_0A_1R_{-1} = \bfzero$ on the right-hand side we deduce that 
$$
P^c \bfx_{-k+m} = \bfzero \Leftrightarrow P \bfx_{-k+m} = \bfx_{-k+m} \Leftrightarrow P \bfx_{-k+(m+1)-1} = \bfx_{-k+(m+1)-1}
$$
which means our supposition is also true for $n = m+1$.  Since the supposition is true for $n=1$ it follows that it is true for each $n=1,\ldots,k$.  Setting $n=k$ shows that $\bfx_{-1} \in U_{\mbox{\scriptsize sg},\,k} \Rightarrow P \bfx_{-1} = \bfx_{-1}$.  Hence $\bfx_{-1} \in P(X)$.  Thus we can state our first new result.

\begin{lemma}
\label{lem2}
{\rm  For all $k \in {\mathbb N}$ we have $U_{\mbox{\scriptsize sg},\,k} \subseteq P(X)$. $\hfill \Box$}
\end{lemma}

\subsection{The finite-length regular Jordan chains}
\label{flrjc}

Once again we follow the methodology of Bart and Lay \cite{bar1} but we reverse the roles of $A_0$ and $A_1$.   Thus we define subspaces $U_{\mbox{\scriptsize rg},\,0} = \{ \bfzero \} \subseteq X$ and
$$
U_{\mbox{\scriptsize rg},\,k} = A_1^{-1}(A_0A_1^{-1})^{k-1}( \{ \bfzero \}) \subseteq X
$$
for each $k \in {\mathbb N}$.  We also define corresponding subspaces $V_{\mbox{\scriptsize rg},\,0} = \{ \bfzero \} \subseteq Y$ and
$$
V_{\mbox{\scriptsize rg},\,k} = (A_0 A_1^{-1})^k( \{ \bfzero \}) \subseteq Y
$$
for each $k \in {\mathbb N}$.  We will use these subspaces to define the finite-length regular Jordan chains of vectors for the pencil $A(z) = A_0 + A_1z$.  Now we observe that $\bfx_1 \in U_{\mbox{\scriptsize rg},\,k} \Leftrightarrow A_1 \bfx_1 \in V_{\mbox{\scriptsize rg},\,k-1} = A_0(U_{\mbox{\scriptsize rg},\,k-1})$.  Thus there exists $\bfx_2 \in U_{\mbox{\scriptsize rg},\,k-1}$ such that $A_1 \bfx_1 = A_0 (- \bfx_2) \Leftrightarrow A_1 \bfx_1 + A_0 \bfx_2 = \bfzero$.  By repeating this argument it follows that we can find a Jordan chain $\{ \bfx_n \}_{n=1}^{k}$ with $\bfx_n \in U_{\mbox{\scriptsize rg},\,k-n+1}$ for each $n=1,\ldots,k$ such that
\begin{equation}
\label{enr}
A_1 \bfx_n + A_0 \bfx_{n+1} = \bfzero
\end{equation}
for each $n = 1,\ldots,k$.  Since $\bfx_k \in U_{\mbox{\scriptsize rg},\,1} = A_1^{-1}( \{ \bfzero\})$ it follows that
\begin{equation}
\label{ekr}
A_1 \bfx_k = \bfzero.
\end{equation}
We say that $\{ \bfx_n \}_{n=1}^{k}$ is a {\em regular} Jordan chain of length $k$ for $A(z)$ and that the chain is generated by the point $\bfx_1 \in U_{\mbox{\scriptsize rg},\,k}$.  The subspace $U_{\mbox{\scriptsize rg},\,k}$ is the generating subspace for all regular Jordan chains of length $k$ for the pencil $A(z)$.  We say these chains are regular because the generating sets are closely related to the coefficients of the nonnegative powers in the Laurent series expansion for $R(z) = A(z)^{-1}$ on the region $z \in {\mathcal U}_{\, s,r}$.  Note that $A_1^{-1}( \{ \bfzero\}) = U_{\mbox{\scriptsize rg},\,1} \subseteq U_{\mbox{\scriptsize rg},\,2} \subseteq \cdots \subseteq U_{\mbox{\scriptsize rg},\,k} \subseteq X$.   We have established the following result.

\begin{lemma}
\label{lem3}
{\rm  For each $k \in {\mathbb N}$ the subspace $U_{\mbox{\scriptsize rg},\,k} = A_1^{-1}(A_0A_1^{-1})^{k-1}( \{ \bfzero \}) \subseteq X$ is the generating subspace for all regular Jordan chains of length $k$ for the pencil $A(z)$.  Thus we have 
$$
U_{\mbox{\scriptsize rg},\,k} = \{ \bfx_1 \in X \mid\ \mbox{there exists}\ \{ \bfx_n \}_{n=1}^k\ \mbox{with}\ \bfx_n \in U_{\mbox{\scriptsize rg},\,k-n+1}\ \mbox{such that (\ref{enr}) and (\ref{ekr}) are satisfied} \} 
$$
for each $k \in {\mathbb N}$. $\hfill \Box$}
\end{lemma}

Consider a regular Jordan chain $\{\bfx_n\}_{n=1}^k$ of length $k$ generated by a point $\bfx_1 \in U_{\mbox{\scriptsize rg},\,k}$.  Now
$$
\bfx_k \in U_{\mbox{\scriptsize rg},\,1} \Rightarrow A_1 \bfx_k = \bfzero \Rightarrow R_{-1}A_1 \bfx_k = \bfzero \Leftrightarrow P\bfx_k = \bfzero \Leftrightarrow P^c \bfx_k = \bfx_k.
$$
Suppose $P^c\bfx_{k-n+1} = \bfx_{k-n+1} \Leftrightarrow R_0A_0 \bfx_{k-n+1} = \bfx_{k-n+1}$ for some $n = m \in {\mathbb N}$ with $1 \leq m < k$.  The sequence $\{\bfx_n\}_{n=1}^k$ is a regular Jordan chain and so $A_1 \bfx_{k-m} + A_0 \bfx_{k-m+1} = \bfzero$.  Therefore $R_0A_1 \bfx_{k-m} + R_0A_0 \bfx_{k-m+1} = \bfzero$ and our supposition that $R_0A_0 \bfx_{k-m+1} = \bfx_{k-m+1}$ now means that
$$
\bfx_{k-m+1} = (-1) R_0A_1 \bfx_{k-m}.
$$
The relationship $A_1 \bfx_{k-m} + A_0 \bfx_{k-m+1} = \bfzero$ also implies $R_{-1}A_1 \bfx_{k-m} + R_{-1} A_0 \bfx_{k-m+1} = \bfzero$ which in turn gives
$$
R_{-1} A_1 \bfx_{k-m} = (-1) R_{-1}A_0 \bfx_{k-m+1} = (-1)^2 R_{-1}A_0R_0A_0 \bfx_{k-m+1}.
$$
By applying the identity $R_{-1}A_0R_0 = \bfzero$ on the right-hand side we deduce that 
$$
P \bfx_{k-m} = \bfzero \Leftrightarrow (I-P) \bfx_{k-m} = \bfx_{k-m} \Leftrightarrow P^c \bfx_{k-(m+1)+1} = \bfx_{k-(m+1)+1}
$$
which means our supposition is also true for $n = m+1$.  Since the supposition is true for $n=1$ it follows that it is true for each $n=1,\ldots,k$.  Setting $n=k$ shows that $\bfx_1 \in U_{\mbox{\scriptsize rg},\,k} \Rightarrow P^c \bfx_1 = \bfx_1$.  Hence $\bfx_1 \in P^c(X)$.  Thus we can state our second new result.

\begin{lemma}
\label{lem4}
{\rm  For all $k \in {\mathbb N}$ we have $U_{\mbox{\scriptsize rg},\,k} \subseteq P^c(X)$. $\hfill \Box$}
\end{lemma}

\subsection{The infinite-length singular Jordan chains}
\label{ilsjc}

When $R(z) = A(z)^{-1}$ has a pole of order $k$ at $z = 0$ and $\{ \bfx_{-n}\}_{n \in {\mathbb N}}$ is a singular Jordan chain generated by a point $\bfx_{-1} \in U_{\mbox{\scriptsize sg},\,k}$ then $\bfx_{-n} = \bfzero$ for all $n \in {\mathbb N}$ with $n > k$.  Consequently we say that the chain has length $k$.  This means that when $R(z)$ has an essential singularity at $z = 0$ we might reasonably expect to find singular Jordan chains of every finite length or indeed of infinite length.  Definitions \ref{def1} and \ref{def2} are new.

\begin{definition}
\label{def1}
{\rm Let $0 \leq s < r \leq \infty$ and suppose $R(z)$ is analytic for all $z \in {\mathcal U}_{s,r}$.  If $\{ \bfx_{-n} \}_{n \in {\mathbb N}} \subseteq X$ with $A_0 \bfx_{-n} + A_1 \bfx_{-n-1} = \bfzero$ for all $n \in {\mathbb N}$ and $\lim_{n \rightarrow \infty} \| \bfx_{-n} \|^{1/n} = s$ then we say that $\{ \bfx_{-n}\}_{n \in {\mathbb N}}$ is an infinite-length singular Jordan chain for the pencil $A(z) = A_0 + A_1z$ on the annular region ${\mathcal U}_{s,\infty}$.  Note that for all $k \in {\mathbb N}$ and $s > 0$ this definition excludes all singular Jordan chains $\{ \bfx_{-n}\}_{n \in {\mathbb N}}$ of length $k$ with $\bfx_{-n} = \bfzero$ for all $n > k$.  When $s = 0$ the finite-length singular Jordan chains are included. } $\hfill \Box$
\end{definition}

For singular chains of length $k$ we started with $U_{\mbox{\scriptsize sg},\,0} = \{ \bfzero\}$ and constructed an increasing sequence $\{U_{\mbox{\scriptsize sg},\,n}\}_{n=1}^k$ of subspaces to find the generating subspace $U_{\mbox{\scriptsize sg},\,k} = A_0^{-1}(A_1A_0^{-1})^{k-1}( \{ \bfzero \}) \subseteq X$.  Although it may seem tempting to define the generating subspace for the infinite-length singular Jordan chains as the subspace $U_{\mbox{\scriptsize sg},\,\infty} = \lim_{k \rightarrow \infty} U_{\mbox{\scriptsize sg},\,k} = \bigcup_{k \in {\mathbb N}} U_{\mbox{\scriptsize sg},\,k}$ it is also necessary to incorporate the limiting conditions on the magnitude of the vectors in the chain and so we prefer a more direct approach.

\begin{definition}
\label{def2}
{\rm  Let $0 \leq s < r \leq \infty$ and suppose $R(z)$ is analytic for all $z \in {\mathcal U}_{s,r}$.   Let $X_s \subseteq X$ be the subspace defined by $X_s = \{ \bfx_{-1} \in X \mid \mbox{there exists}\ \{ \bfx_{-n} \}_{n \in {\mathbb N}}\ \mbox{with}\ A_0 \bfx_{-n} + A_1 \bfx_{-n-1} = \bfzero\ \mbox{for all}\ n \in {\mathbb N}\ \mbox{and}\ \lim_{n \rightarrow \infty} \| \bfx_{-n} \|^{1/n} = s \}$.  We say that $X_s$ is the generating set for all infinite-length singular Jordan chains on the region ${\mathcal U}_{s,\infty}$.}  $\hfill \Box$
\end{definition} 

Consider an infinite-length Jordan chain $\{\bfx_{-n}\}_{n \in {\mathbb N}}$ generated by a point $\bfx_{-1} \in X$.  Since $A_0 \bfx_{-1} + A_1 \bfx_{-2} = \bfzero$ we have $R_{-1}A_0 \bfx_{-1} + R_{-1}A_1 \bfx_{-2} = \bfzero$ and we can use the identity $R_{-1}A_0R_0 = 0$ to deduce that
$$
R_{-1}A_1 \bfx_{-2} =  (-1) R_{-1}A_0 \bfx_{-1} =  (-1) (R_{-1}A_0) \left[ R_0A_0 + R_{-1}A_1\right] \bfx_{-1} = (-1) (R_{-1}A_0) R_{-1}A_1 \bfx_{-1}.
$$
Now suppose for some $k \in {\mathbb N}$ we have
$$
R_{-1}A_1 \bfx_{-k} = (-1)^{k-1}(R_{-1}A_0)^{k-1} R_{-1}A_1 \bfx_{-1}.
$$
Since $A_0 \bfx_{-k} + A_1 \bfx_{-k-1} = \bfzero$ we have $R_{-1}A_0 \bfx_{-k} + R_{-1}A_1 \bfx_{-k-1} = \bfzero$ and so
\begin{eqnarray*}
R_{-1}A_1 \bfx_{-k-1} =  (-1) R_{-1}A_0 \bfx_{-k} & = & (-1) (R_{-1}A_0) \left[ R_0A_0 + R_{-1}A_1 \right] \bfx_{-k} \\
& = & (-1) (R_{-1}A_0) R_{-1}A_1 \bfx_{-k} = (-1)^k (R_{-1}A_0)^k R_{-1}A_1 \bfx_{-1}
\end{eqnarray*}
where we have once again used the identity $R_{-1}A_0R_0 = 0$.  It follows by induction that $R_{-1}A_1\bfx_{-n} = (-1)^{n-1}(R_{-1}A_0)^{n-1} R_{-1}A_1\bfx_{-1}$ for all $n \in {\mathbb N}$.

Now consider an alternative viewpoint.  Choose $n \in {\mathbb N}$ with $n > 1$.  Since $A_0 \bfx_{-n+1} + A_1 \bfx_{-n} = \bfzero$ we have $R_0A_0 \bfx_{-n+1} + R_0A_1 \bfx_{-n} = \bfzero$ and hence we can use the identity $R_0A_1R_{-1} = 0$ to show that
$$
R_0A_0 \bfx_{-n+1} = (-1) R_0A_1 \bfx_{-n} = (-1) (R_0A_1) \left[ R_0A_0 + R_{-1}A_1\right] \bfx_{-n} = (-1) (R_0A_1) R_0A_0 \bfx_{-n}.
$$
Now suppose for some $k \in {\mathbb N}$ with $k < n$ we have
$$
R_0A_0 \bfx_{-n+k-1} = (-1)^{k-1}(R_0A_1)^{k-1} R_0A_0 \bfx_{-n}.
$$
Since $A_0 \bfx_{-n+k} + A_1 \bfx_{-n+k-1} = \bfzero$ we have $R_0A_0 \bfx_{-n+k} + R_0A_1 \bfx_{-n+k-1} = \bfzero$ and we can see that
\begin{eqnarray*}
R_0A_0 \bfx_{-n+k} = (-1) R_0A_1 \bfx_{-n+k-1} & = & (-1) (R_0A_1) \left[ R_0A_0 + R_{-1}A_1 \right] \bfx_{-n+k-1} \\
& = & (-1) (R_0A_1) R_0A_0 \bfx_{-n+k-1} = (-1)^k (R_0A_1)^k R_0A_0\bfx_{-n}
\end{eqnarray*}
where we have once again used the identity $R_0A_1R_{-1} = 0$.  It follows by induction that $R_0A_0\bfx_{-n+k-1} = (-1)^{k-1}(R_0A_1)^{k-1} R_0A_0 \bfx_{-n}$ for all $k =1,\ldots,n$.  

In particular we have established that if $\bfx_{-1} \in X$ and $\{\bfx_{-n}\}_{n \in {\mathbb N}}$ is an associated infinite-length Jordan chain then for each $n \in {\mathbb N}$ we have
\begin{equation}
\label{jcf1}
R_{-1}A_1\bfx_{-n} = (-1)^{n-1}(R_{-1}A_0)^{n-1} R_{-1}A_1\bfx_{-1} \Leftrightarrow P \bfx_{-n} = R_{-n} A_1 P \bfx_{-1}
\end{equation}
where we have used $R_{-1}A_1 + R_0A_0 = I$ and $R_{-1}A_1R_0 = 0$ to justify the equivalence.  We also have
\begin{equation}
\label{jcf2}
R_0A_0\bfx_{-1} = (-1)^{n-1}(R_0A_1)^{n-1} R_0A_0 \bfx_{-n} \Leftrightarrow P^c\bfx_{-1} = R_{n-1}A_0 P^c\bfx_{-n}.
\end{equation}
 Now suppose that $\bfx_{-1} \in X_s$ and assume that $\|P^c \bfx_{-1} \| = c > 0$.  Since $\lim_{n \rightarrow \infty} \|\bfx_{-n}\|^{1/n} = s$ it follows that for each $\delta > 0$ there exists a real constant $c_{\delta} > 0$ such that $\| \bfx_{-n}\| < c_{\delta} (s + \delta)^n$ for all $n \in {\mathbb N}$.  Also, since $\lim_{n \rightarrow \infty} \|R_{-n} A_1 \|^{1/n} = s$, there is a real constant $d_{\delta} > 0$ such that $\|P \bfx_{-n}\| = \|R_{-n}A_1 P \bfx_{-1} \| < d_{\delta} (s+\delta)^n$ for all $n \in {\mathbb N}$.  Therefore
$$
\| P^c \bfx_{-n} \| \leq \| \bfx_{-n} \| + \| P \bfx_{-n} \| \leq (c_{\delta} + d_{\delta}) (s + \delta)^n
$$
for all $n \in {\mathbb N}$.  Since $\lim_{n \rightarrow \infty} \|R_{n-1}A_0\|^{1/n} = 1/r$ it follows from (\ref{jcf2}) that for each $\epsilon > 0$ and sufficiently small we can find a real constant $f_{\epsilon} > 0$ such that
$$
c = \|P^c \bfx_{-1} \| = \|R_{n-1}A_0 P^c\bfx_{-n}\| \leq f_{\epsilon} \left[1 / (r - \epsilon)^n \right] \|P^c \bfx_{-n} \|
$$
for all $n \in {\mathbb N}$.  Combining this inequality with the previous inequality gives
$$
\| P^c \bfx_{-n} \| \leq (1/c) (c_{\delta} + d_{\delta}) f_{\epsilon} \left[ (s + \delta) / (r - \epsilon) \right]^n \|P^c \bfx_{-n} \|
$$
for all $n \in {\mathbb N}$.  If we choose $\delta$ and $\epsilon$ small enough to ensure that $s + \delta < r - \epsilon$ and if we subsequently choose $n \in {\mathbb N}$ so that $(1/c) (c_{\delta} + d_{\delta}) f_{\epsilon} \left[ (s + \delta) / (r - \epsilon) \right]^n < 1/2$ then it follows that $\|P^c \bfx_{-n} \| \leq (1/2) \|P^c \bfx_{-n} \|$.  Therefore $\|P^c \bfx_{-n} \| = 0 \Leftrightarrow P^c \bfx_{-n} = \bfzero \Leftrightarrow P \bfx_{-n} = \bfx_{-n} \Rightarrow P \bfx_{-1} = \bfx_{-1} \iff P^c \bfx_{-1} = \bfzero$.  This contradicts our assumption that $\|P^c \bfx_{-1}\| = c > 0$.   Therefore $P^c \bfx_{-1} = \bfzero \iff P \bfx_{-1} = \bfx_{-1}$.  Thus we have shown that $X_s \subseteq P(X)$.  The next two results are new.

\begin{lemma}
\label{lem5}
{\rm  Let $0 \leq s < r \leq \infty$ and suppose $R(z)$ is analytic on ${\mathcal U}_{s,\, r}$.  Then $X_s = P(X)$. $\hfill \Box$  }
\end{lemma}

{\bf Proof.}  We already showed that $X_s \subseteq P(X)$.   Now suppose $\bfx_{-1} \in X_s^a$ where we use the notation $E^a \subseteq X$ for each subset $E \subseteq X$ to denote the closure of $E$ in $X$.   It follows that for each $\epsilon > 0$ we can find a point $\bfx_{-1,\epsilon} \in X_s$ with $\| \bfx_{-1,\epsilon} - \bfx_{-1} \| < \epsilon$ such that $\bfx_{-1,\epsilon}$ generates an infinite-length singular Jordan chain for $A(z)$ on ${\mathcal U}_{s,\infty}$.  Therefore $P \bfx_{-1,\epsilon} = \bfx_{-1,\epsilon}$.  Since $\|P\| = \|R_{-1}A_1\| < \infty$ it follows that
$$
\| P \bfx_{-1} - \bfx_{-1} \| = \| P(\bfx_{-1} - \bfx_{-1,\epsilon}) + (\bfx_{-1,\epsilon} - \bfx_{-1}) \| \leq (\|R_{-1}A_1\| + 1) \| \bfx_{-1,\epsilon} - \bfx_{-1}\| \leq (\|P\| + 1) \epsilon
$$
and since $\epsilon > 0$ is arbitrary it follows that $\| P \bfx_{-1} - \bfx_{-1} \| = 0$.  Hence $P \bfx_{-1} = \bfx_{-1}$.  Therefore $\bfx_{-1} \in P(X)$.  Thus we also have $X_s^a \subseteq P(X)$.  Conversely, if $\bfx_{-1} \in P(X)$ and we define $\bfx_{-n} = R_{-n} \bfx_{-1}$ for all $n \in {\mathbb N}$, the identities $R_{-n} = (-1)^{n-1}(R_{-1}A_0)^{n-1}R_{-1}$ and $R_0A_0R_{-1} = 0$ established in \cite{alb2} ensure that $P^c \bfx_{-n} = \bfzero$ for all $n \in {\mathbb N}$.  It follows from the formula for $\bfx_{-n}$ that
$$
A_0 \bfx_{-n} + A_1 \bfx_{-n-1} = A_0 \bfx_{-n+1} - A_1 (R_{-1}A_0) \bfx_{-n} = (I - A_1R_{-1})A_0 \bfx_{-n} = A_0 R_0 A_0 \bfx_{-n} = A_0P^c \bfx_{-n} = \bfzero
$$
for all $n \in {\mathbb N}$.  Since $\bfx_{-n} = R_{-n} \bfx_{-1}$ and $\lim_{n \rightarrow \infty} \|R_{-n}\|^{1/n} = s$ we also have $\lim_{n \rightarrow \infty} \| \bfx_{-n} \|^{1/n} = s$.  Thus $\{ \bfx_{-n} \}_{n \in {\mathbb N}}$ is an infinite-length singular Jordan chain for $A(z)$ on ${\mathcal U}_{s,\infty}$.   Therefore $\bfx_{-1} \in X_s$.  Hence $P(X) \subseteq X_s$.

Thus we have shown that $X_s \subseteq X_s^a \subseteq P(X) \subseteq X_s$.  It follows that $X_s = X_s^a = P(X)$.  $\hfill \Box$

\begin{corollary}
\label{cor1}
{\rm  The set $X_s = X_s^{a}$ is closed.} $\hfill \Box$
\end{corollary}

\subsection{The infinite-length regular Jordan chains}
\label{ilrjc}

We wish to obtain a similar characterization of the subspace $P^c(X)$.  Definitions \ref{def3} and \ref{def4} are new.

\begin{definition}
\label{def3}
{\rm Let $0 \leq s < r \leq \infty$ and suppose $R(z)$ is analytic for all $z \in {\mathcal U}_{s,r}$.  If $\{ \bfx_{n} \}_{n \in {\mathbb N}} \subseteq X$ with $A_0 \bfx_{n} + A_1 \bfx_{n-1} = \bfzero$ for all $n \in {\mathbb N}$ and $\lim_{n \rightarrow \infty} \| \bfx_{n} \|^{1/n} = 1/r$ then we say that $\{ \bfx_{n}\}_{n \in {\mathbb N}}$ is an infinite-length regular Jordan chain for the pencil $A(z) = A_0 + A_1z$ on the annular region ${\mathcal U}_{0,r}$.  Note that for all $k \in {\mathbb N}$ and $r < \infty$ this definition excludes all regular Jordan chains $\{ \bfx_{n}\}_{n \in {\mathbb N}}$ of length $k$ with $\bfx_{n} = \bfzero$ for all $n > k$.  When $r = \infty$ the finite-length regular Jordan chains are included.} $\hfill \Box$
\end{definition}

\begin{definition}
\label{def4}
{\rm  Let $0 \leq s < r \leq \infty$ and suppose $R(z)$ is analytic for all $z \in {\mathcal U}_{s,r}$.  Let $X_r \subseteq X$ be the subspace defined by $X_r = \{ \bfx_1 \in X \mid \mbox{there exists}\ \{ \bfx_n \}_{n \in {\mathbb N}}\ \mbox{with}\ A_1 \bfx_n + A_0 \bfx_{n+1} = \bfzero\ \mbox{for all}\ n \in {\mathbb N}\ \mbox{and}\  \lim_{n \rightarrow \infty} \| \bfx_n \|^{1/n} = 1/r \}$.  We say that $X_r$ is the generating set for all infinite-length regular Jordan chains on the region ${\mathcal U}_{0,r}$.}  $\hfill \Box$
\end{definition}
  
We can now use some elementary algebra to show that regular Jordan chains are structurally the same as singular Jordan chains.  Let $B(w) = B_0 + B_1w$ where $B_0 = A_1$ and $B_1 = A_0$ and $w = 1/z \in {\mathbb C}$.  It follows that $w \in {\mathcal U}_{1/r,1/s}$ if and only if $z \in {\mathcal U}_{s,r}$.  If we define $S(w) = B(w)^{-1}$ then $S(w) = zR(z)$ and if we write $S(w) = \sum_{k \in {\mathbb Z}} S_k w^k$ for each $w \in {\mathcal U}_{1/r,1/s}$ then $S_k = R_{-k-1}$ for all $k \in {\mathbb Z}$.  In particular $P = R_{-1}A_1 = S_0B_0$ and $P^c = R_0A_0 = S_{-1}B_1$.  If we set $\bfu_{-1} = \bfx_1 \in X_r$ and define $\bfu_{-n} = \bfx_n$ for all $n \in {\mathbb N}$ where $\{\bfx_n\}_{n \in {\mathbb N}}$ is the infinite-length regular Jordan chain for $A(z)$ generated by $X_r$ with $\lim_{n \rightarrow \infty} \| \bfx_n \|^{1/n} = 1/r$ then we can see that $\{\bfu_{-n}\}_{n \in {\mathbb N}} = \{\bfx_n\}_{n \in {\mathbb N}}$ is an infinite-length singular Jordan chain for $B(w)$ generated by $X_r$ with $\lim_{n \rightarrow \infty} \| \bfu_{-n}\|^{1/n} = 1/r$.  Now we can establish two more new results.

\begin{lemma}
\label{lem6}
{\rm  Let $0 \leq s < r \leq \infty$ and suppose $R(z)$ is analytic on ${\mathcal U}_{s,\, r}$.  Then $X_r = P^c(X)$}. $\hfill \Box$
\end{lemma}

{\bf Proof.}  Apply the same arguments used to establish Lemma \ref{lem5} with $A(z)$ replaced by $B(w)$, $R_0$ replaced by $S_0$, $R_{-1}$ replaced by $S_{-1}$, $X_s$ replaced by $X_r$ and $\{ \bfx_{-n}\}_{n \in {\mathbb N}}$ replaced by $\{ \bfu_{-n}\}_{n \in {\mathbb N}}$.  Note that $P = R_{-1}A_1$ is replaced by $P^c = S_{-1}B_1$.  $\hfill \Box$ 

\begin{corollary}
\label{cor2}
{\rm  The set $X_r = X_r^{a}$ is closed.} $\hfill \Box$
\end{corollary}

\subsection{The domain space decomposition}
\label{dsd}

Let $0 \leq s < r \leq \infty$ and suppose that $R(z)$ is analytic for $z \in {\mathcal U}_{s,r}$.  Let $P = R_{-1}A_1 \in {\mathcal B}(X)$ and $P^c = R_0A_0 \in {\mathcal B}(X)$ be the key projections that separate the bounded and unbounded parts of the spectral set for $A(z)$ relative to the annular region ${\mathcal U}_{s,r}$ and let $X_s$ and $X_r$ be the respective generating subspaces for the corresponding infinite-length singular and regular Jordan chains for $A(z)$ on ${\mathcal U}_{\, s,r}$.  The results in Sections \ref{flsjc}, \ref{flrjc}, \ref{ilsjc} and \ref{ilrjc} show that $X_s = P(X)$ and $X_r = P^c(X) = (I - P)(X)$ and hence that $X = X_s \oplus X_r \cong X_s \times X_r$.  This is the first part of our Theorem \ref{th1}.

\subsection{The range space decomposition}
\label{rsd}

Thus far we have concentrated on developing a structural decomposition for the domain space $X$.  We will now show that the range space $Y$ has a corresponding structural decomposition.  Define $Y_s = Q(Y)$ and $Y_r = Q^c(Y)$.  Therefore $Y = Y_s \oplus Y_r \cong Y_s \times Y_r$.  Since we have assumed that $R(z) = A(z)^{-1}$ is well defined for $z \in {\mathcal U}_{\, s,r}$ we know that $A(z) = (A_0 + A_1z)$ is a $1$\textendash $1$ mapping of $X$ onto $Y$ for all $z \in {\mathcal U}_{\, s,r}$.  While it is obvious that $QA_1 = A_1R_{-1}A_1 = A_1P$ and $Q^cA_0 = A_0R_0A_0 = A_0P^c$ it is also true that $QA_0 = (I -A_0R_0)A_0 = A_0(I - R_0A_0) = A_0P$ and $Q^cA_1 = (I -A_1R_{-1})A_1 = A_1(I - R_{-1}A_1) = A_1P^c$.   Therefore
$$
Y_s =  Q(Y) = Q(A_0 + A_1z)(X) = (A_0 + A_1z)P(X) = (A_0 + A_1z)X_s
$$
and
$$
Y_r = Q^c(Y) = Q^c(A_0 + A_1z)(X) = (A_0 + A_1z)P^c(X) = (A_0 + A_1z)X_r.
$$
Hence $(A_0+A_1z)$ is a $1$\textendash $1$ mapping of $X_s$ onto $Y_s$ and a $1$\textendash $1$ mapping of $X_r$ onto $Y_r$ for all $z \in {\mathcal U}_{\, s,r}$.  More specifically we can also see that $Q \bfy = A_1(P\bfx)$ if and only if $R_{-1}(Q\bfy) = P \bfx$ and that $Q^c \bfy = A_0(P^c \bfx)$ if and only if $R_0(Q^c \bfy) = P^c \bfx$.  Thus $A_1[P(X)] = Q(Y)$, $R_{-1}[Q(Y)] = P(X)$, $A_0[P^c(X)] = Q^c(Y)$ and $R_0[Q^c(Y)] = P^c(X)$.   It follows that $A_1$ is a $1$\textendash $1$ mapping of $X_s$ onto $Y_s$ and $A_0$ is a $1$\textendash $1$ mapping of $X_r$ onto $Y_r$ while $R_{-1}$ is a $1$\textendash $1$ mapping of $Y_s$ onto $X_s$ and $R_0$ is a $1$\textendash $1$ mapping of $Y_r$ onto $X_r$.  Thus we have shown that $Y_s = Q(Y) = A_1(X_s)$ and $Y_r = Q^c(Y) = (I - Q)(Y) = A_0(X_r)$ and hence that $Y = Y_s \oplus Y_r \cong Y_s \times Y_r$.  This is the second part of our Theorem \ref{th1}.

\subsection{Separation of the fundamental equations}
\label{sfe}

The projections $P=R_{-1}A_1$ and $Q=A_1R_{-1}$ separate each of the doubly-infinite left and right systems of fundamental equations into two separate semi-infinite left and right systems.  Each of these semi-infinite systems can be solved recursively.  See \cite{alb3} for more details.  In fact if we define ${\mathfrak A}_i = QA_iP \in {\mathcal B}(X_s,Y_s)$ and ${\mathfrak A}_i^c = Q^cA_iP^c \in {\mathcal B}(X_r,Y_r)$ for each $i=0,1$ and ${\mathfrak R}_j = PR_jQ \in {\mathcal B}(Y_r,X_r)$ and ${\mathfrak R}_j^c = P^cR_jQ^c \in {\mathcal B}(Y_s,X_s)$ for each $j=-1,0$ then ${\mathfrak R}_{-1}^c = 0$ and ${\mathfrak R}_0 = 0$ and the identities that define the basic solution $\{R_{-1},R_0\} \cong \{{\mathfrak R}_{-1}, {\mathfrak R}_0^c\}$ take the simplified form
\begin{description}
\item[{\bf (a)}] ${\mathfrak R}_{-1}{\mathfrak A}_1 = {\mathfrak J}$ and ${\mathfrak R}_0^c{\mathfrak A}_0^c = {\mathfrak J}^c$ for the left system; and
\item[{\bf (b)}] ${\mathfrak A}_1{\mathfrak R}_{-1} = {\mathfrak L}$ and ${\mathfrak A}_0^c{\mathfrak R}_0^c = {\mathfrak L}^c$ for the right system;
\end{description}
where ${\mathfrak J} \in {\mathcal B}(X_s)$, ${\mathfrak J}^c \in {\mathcal B}(X_r)$, ${\mathfrak L} \in {\mathcal B}(Y_s)$ and ${\mathfrak L}^c \in {\mathcal B}(Y_r)$ are the relevant identity mappings.

\section{The solution process}
\label{sp}

Let us now consider the solution process.  In this section we will assume that $X$ has a Schauder basis $\{ \bfe_n\}_{n \in {\mathbb N}}$ where we assume without loss of generality that $\|\bfe_n\| = 1$ for all $n \in {\mathbb N}$.  Define $\bfp_n = P\bfe_n$ and $\bfp_n^c = P^c \bfe_n$ for all $n \in {\mathbb N}$.  We have the following result.

\begin{lemma}
\label{lem7}
{\rm  If $\bfp \in P(X)$ then there is a unique sequence $\{ \alpha_n \}_{n \in {\mathbb N}} \subseteq {\mathbb C}$ such that $\bfp = \sum_{n \in {\mathbb N}} \alpha_n \bfp_n$ and $\bfzero = \sum_{n \in {\mathbb N}} \alpha_n \bfp_n^c$.   Likewise if $\bfp^c \in P^c(X)$ then there is a unique sequence $\{ \beta_n \}_{n \in {\mathbb N}} \subseteq {\mathbb C}$ such that $\bfp^c = \sum_{n \in {\mathbb N}} \beta_n \bfp_n^c$ and $\bfzero = \sum_{n \in {\mathbb N}} \beta_n \bfp_n$.} $\hfill \Box$
\end{lemma}

{\bf Proof.} Since $\bfp \in P(X) \subseteq X$ we have
$$
\bfp = \sum_{n \in {\mathbb N}} \alpha_n \bfe_n \Leftrightarrow \sum_{n \in {\mathbb N}} \alpha_n (P\bfe_n + P^c\bfe_n) \Leftrightarrow \sum_{n \in {\mathbb N}} \alpha_n \bfp_n + \sum_{n \in {\mathbb N}} \alpha_n \bfp_n^c \Leftrightarrow \{\bfp = \sum_{n \in {\mathbb N}} \alpha_n \bfp_n\ \mbox{and}\ \bfzero = \sum_{n \in {\mathbb N}} \alpha_n \bfp_n^c\}.
$$ 
A similar argument applies to the representation defined by $\bfp^c = \sum_{n\in {\mathbb N}} \beta_n \bfp_n^c$ and $\bfzero = \sum_{n \in {\mathbb N}} \beta_n \bfp_n$.   $\hfill \Box$

The existence of a unique representation for all vectors in a subspace using a spanning set that is not a basis for the subspace justifies our next definition of a {\em conditional frame}.  Note that this specific term is apparently not used in the established literature \cite{cas1}. 

\begin{definition}
\label{def5}
{\rm If $\{ \bfe_n \}_{n \in {\mathbb N}}$ is a Schauder basis for $X$ and we define $\bfp_n = P\bfe_n$ and $\bfp_n^c = P^c \bfe_n$ for each $n \in {\mathbb N}$ then we say that $\{ \bfp_n\}_{n \in {\mathbb N}} \in P(X)$ defines a conditional frame for $P(X)$ and $\{ \bfp_n^c\}_{n \in {\mathbb N}} \in P^c(X)$ defines a conditional frame for $P^c(X)$.} $\hfill \Box$
\end{definition}

The next result is important because the Schauder basis for $X$ is often constructed from vectors that span the key subspaces.  In fact it is usual to begin the solution process by constructing vectors that generate the desired Jordan chains and hence lie within either $P(X)$ or $P^c(X)$.

\begin{corollary}
\label{cor3}
{\rm Let $\{ \bfe_n \}_{n \in {\mathbb N}}$ be a Schauder basis for $X$ and define $\bfp_n = P \bfe_n$ and $\bfp_n^c = P^c \bfe_n$ for each $n \in {\mathbb N}$.  There exist sets ${\mathbb L} \subseteq {\mathbb N}$ and ${\mathbb M} = {\mathbb N} \setminus {\mathbb L} \subseteq {\mathbb N}$ such that $\{ \bfp_n \}_{n \in {\mathbb L}}$ is a basis for $P(X)$ and $\{ \bfp_n^c \}_{n \in {\mathbb M}}$ is a basis for $P^c(X)$ if and only if $\bfe_n = \bfp_n$ for all $n \in {\mathbb L}$ and $\bfe_n = \bfp_n^c$ for all $n \in {\mathbb M}$. } $ \hfill \Box$
\end{corollary} 

We assume $X$ has a Schauder basis $\{ \bfe_n\}_{n \in {\mathbb N}}$ with $\bfp_n = P\bfe_n$ and $\bfp_n^c = P^c \bfe_n$ for all $n \in {\mathbb N}$.  We recall that $A_1P(X) = Q(Y) \Leftrightarrow P(X) = R_{-1}Q(Y)$ and $A_0P^c(X) = Q^c(Y) \Leftrightarrow P^c(X) = R_0Q^c(Y)$.  Thus we define $\bfq_n = A_1 \bfp_n \Leftrightarrow R_{-1}\bfq_n = \bfp_n$ and $\bfq_n^c = A_0 \bfp_n^c \Leftrightarrow R_0\bfq_n^c = \bfp_n^c$ for each $n \in {\mathbb N}$.  We also define $\bff_n = \bfq_n + \bfq_n^c$ for each $n \in {\mathbb N}$.  We have the following result.

\begin{lemma}
\label{lem8}
{\rm  The set $\{ \bff_n \}_{n \in {\mathbb N}} \subseteq Y$ is a Schauder basis for $Y$.  } $\hfill \Box$
\end{lemma}

{\bf Proof.}  Let $\bfy \in Y$.   Then we can write $\bfy = \bfq + \bfq^c$ where $\bfq \in Q(Y)$ and $\bfq^c \in Q^c(Y)$.  Now define $\bfp = R_{-1}\bfq \Leftrightarrow A_1 \bfp = \bfq$ and choose $\{ \alpha_n \}_{n \in {\mathbb N}} \in {\mathbb C}$ such that $\bfp = \sum_{n \in {\mathbb N}} \alpha_n \bfp_n$ and $\bfzero = \sum_{n \in {\mathbb N}} \alpha_n \bfp_n^c$.  Therefore
$$
\bfq = A_1(\sum_{n \in {\mathbb N}} \alpha_n \bfp_n) = \sum_{n \in {\mathbb N}} \alpha_n A_1 \bfp_n = \sum_{n \in {\mathbb N}} \alpha_n \bfq_n.
$$
We also have 
$$
\bfzero = A_0(\sum_{n \in {\mathbb N}} \alpha_n \bfp_n^c) = \sum_{n \in {\mathbb N}} \alpha_n A_0 \bfp_n^c = \sum_{n \in {\mathbb N}} \alpha_n \bfq_n^c.
$$
A similar argument shows that we can choose $\{ \beta_n \}_{n \in {\mathbb N}} \in {\mathbb C}$ such that $\bfq^c = \sum_{n \in {\mathbb N}} \beta_n \bfq_n^c$ and $\bfzero = \sum_{n \in {\mathbb N}} \beta_n \bfq_n$.  Now it follows that $\bfy = \sum_{n \in {\mathbb N}} \gamma_n \bff_n$ where we have defined $\gamma_n = \alpha_n + \beta_n$ for all $n \in {\mathbb N}$.  To show the representation is unique we assume that
\begin{equation}
\label{gnfnzero}
\bfzero = \sum_{n \in {\mathbb N}} \gamma_n \bff_n.
\end{equation}
If we multiply (\ref{gnfnzero}) on the left by $R_{-1}$ and use the fact that $R_{-1}\bff_n = \bfp_n$ then we obtain $\bfzero = \sum_{n \in {\mathbb N}} \gamma_n \bfp_n$ which means that $\bfx = \sum_{n \in {\mathbb N}} \gamma_n \bfe_n = \sum_{n \in {\mathbb N}} \gamma_n (\bfp_n + \bfp_n^c) = \sum_{n \in {\mathbb N}} \gamma_n \bfp_n^c \in P^c(X)$.  If we multiply (\ref{gnfnzero}) on the left by $R_0$ and use the fact that $R_0 \bff_n = \bfp_n^c$ then we obtain $\bfzero = \sum_{n \in {\mathbb N}} \gamma_n \bfp_n^c$ which means that $\bfx = \sum_{n \in {\mathbb N}} \gamma_n \bfe_n = \sum_{n \in {\mathbb N}} \gamma_n (\bfp_n + \bfp_n^c) = \sum_{n \in {\mathbb N}} \gamma_n \bfp_n \in P(X)$.  Since $\bfx \in P(X)$ and $\bfx \in P^c(X)$ it follows that $\bfx = \bfzero$.  Therefore $\gamma_n = 0$ for all $n \in {\mathbb N}$.   Hence the representation is unique.  $\hfill \Box$.

\begin{corollary}
\label{cor4}
{\rm The set $\{ \bfq_n \}_{n \in {\mathbb N}} \in Q(Y)$ forms a conditional frame for $Q(Y)$ and the set $\{ \bfq_n^c \}_{n \in {\mathbb N}} \in Q^c(Y)$ forms a conditional frame for $Q^c(Y)$. } $\hfill \Box$
\end{corollary} 

In order to find the key projections it is first necessary to find the complete collection of relevant Jordan chains and the associated decompositions.  The basic solution $\{R_{-1},R_0\} \cong \{{\mathfrak R}_{-1}, {\mathfrak R}_0^c \}$ can then be found by solving the equations corresponding to the identities {\bf (a)} and {\bf (b)} described in the final paragraph of Section \ref{ilgjc}.  The Schauder bases and associated conditional frames ensure that these identities are reduced to a semi-infinite set of matrix equations which can, in principle, be solved by Gaussian elimination.

We begin with a Schauder basis $\{\bfe_n\}_{n \in {\mathbb N}}$ for $X$.  Find all singular Jordan chains $\{ \bfx_{-n} \}_{n \in {\mathbb N}}$ on the region ${\mathcal U}_{s,\infty}$ by finding the most general parametric solutions to the system of equations $A_0 \bfx_{-n} + A_1 \bfx_{-n-1} = \bfzero$ for all $n \in {\mathbb N}$ such that $\lim_{n \rightarrow \infty} \| \bfx_{-n}\|^{1/n} = s$.  Let $\{ \bfp_n \}_{n \in {\mathbb L}}$, where ${\mathbb L} \subseteq {\mathbb N}$, be a basis for the space $X_s$ generated by these Jordan chains.  Write
$$
\bfp_n = \sum_{m \in {\mathbb N}} p_{m n}\bfe_m \iff \bfp_n = [ p_{m n} ]_{m \in {\mathbb N}}
$$
for each $n \in {\mathbb L}$ where $[ p_{m n} ]_{m \in {\mathbb N}}$ is the coordinate representation of $\bfp_n$ with respect to the basis $\{ \bfe_n \}_{n \in {\mathbb N}}$.   Find all regular Jordan chains $\{ \bfx_n \}_{n \in {\mathbb N}}$ on the region ${\mathcal U}_{0,r}$ by finding the most general parametric solutions to the system of equations $A_0 \bfx_n + A_1 \bfx_{n-1} = \bfzero$ for all $n \in {\mathbb N}-1$ such that $\lim_{n \rightarrow \infty} \| \bfx_n \|^{1/n} = 1/r$.  Let $\{ \bfp_n^c \}_{n \in {\mathbb M}}$ where ${\mathbb M} \subseteq {\mathbb N}$ be a basis for the space $X_r$ generated by these Jordan chains.  Write
$$
\bfp_n^c = \sum_{m \in {\mathbb N}} p_{mn}^c \bfe_m \iff \bfp_n^c = [ p_{m n}^c ]_{m \in {\mathbb N}}
$$
for each $n \in {\mathbb M}$ where $[ p_{m n}^c ]_{m \in {\mathbb N}}$ is the coordinate representation of $\bfp_n^c$ with respect to the basis $\{ \bfe_n \}_{n \in {\mathbb N}}$.   Since $X = X_s \oplus X_r$, Corollary \ref{cor3} shows that ${\mathbb M} = {\mathbb N} \setminus {\mathbb L}$.  We will now write
$$
\bxi = \sum_{n \in {\mathbb L}} \xi_n \bfp_n \iff \bxi = [ \xi_n ]_{n \in {\mathbb L}} \quad \mbox{and} \quad \bfeta = \sum_{n \in {\mathbb M}} \eta_n \bfp_n^c \iff \bfeta = [ \eta_n ]_{n \in {\mathbb M}}
$$
where $[ \xi_n ]_{n \in {\mathbb L}}$ is the coordinate representation of $\bxi$ relative to the basis $\{\bfp_n\}_{n \in {\mathbb L}}$ and $[ \eta_n ]_{n \in {\mathbb M}}$ is the coordinate representation of $\bfeta$ relative to the basis $\{\bfp_n^c\}_{n \in {\mathbb M}}$.  We can now define linear transformations $L = [\bfp_1,\bfp_2,\ldots]  \in {\mathcal B}(X_s, X)$ by setting $L \bxi = \sum_{n \in {\mathbb L}} \xi_n \bfp_n$ for each $\bxi \in X_s$ and  $M = [ \bfp_1^c, \bfp_2^c, \ldots ] \in {\mathcal B}(X_r,X)$ by setting $M \bfeta = \sum_{n \in {\mathbb M}} \eta_n \bfp_n^c$ for each $\bfeta \in X_r$ and an augmented linear transformation $E = [L \mid M] \in {\mathcal B}(X_s \times X_r,X)$ defined by
$$
E \blambda = [L \mid M] \left[ \begin{array}{c}
\bxi \\ \hline
\bfeta \end{array} \right] = L\bxi + M \bfeta
$$ 
for each $\bxi \in X_s$ and $\bfeta \in X_r$.  Since $\{ \bfp_n \}_{n \in {\mathbb L}} \subseteq X_s$ and $\{ \bfp_n^c \}_{n \in {\mathbb M}} \subseteq X_r$ are basis sets for the respective subspaces $X_s$ and $X_r$ and since $X = X_s \oplus X_r  \cong X_s \times X_r$ it follows that the vectors $\{ \{ \bfp_n \}_{n \in {\mathbb L}}, \{ \bfp_n^c \}_{n \in {\mathbb M}} \}$ form a basis for $X$.  Therefore the transformation $E \in {\mathcal B}(X_s \times X_r,X)$ is invertible.  If we write
$$
E^{-1} = \left[ \begin{array}{c}
F \\ \hline
G \end{array} \right] \in {\mathcal B}(X, X_s \times X_r)
$$
then it follows that $LF + MG = I$ and
$$
\left[ \begin{array}{c|c}
FL & FM \\ \hline
GL & GM \end{array} \right] = \left[ \begin{array}{c|c}
I & 0 \\ \hline
0 & I \end{array} \right].
$$
For each $ \bfx \in X$ we wish to find $\blambda \in X_s \times X_r$ such that
$$
\bfx = E \blambda =  [L \mid M] \left[ \begin{array}{c}
\bxi \\ \hline
\bfeta \end{array} \right] = L \bxi + M \bfeta \iff \blambda = E^{-1} \bfx \iff \left[ \begin{array}{c}
\bxi \\ \hline
\bfeta \end{array} \right] = \left[ \begin{array}{c}
F\bfx \\ \hline
G\bfx \end{array} \right].
$$
Thus we have $\bfx = LF \bfx + MG \bfx$ where $P = LF \in {\mathcal B}(X,X_s)$ and $P^c = MG \in {\mathcal B}(X,X_r)$ are the desired key projections.  In the image space we can define $\{\bfq_n\}_{n \in {\mathbb L}} \in Y_s$ by setting $\bfq_n = A_1 \bfp_n$ and $\{ \bfq_n^c\}_{n \in {\mathbb M}} \in Y_r$ by setting $\bfq_n^c = A_0 \bfp_n^c$.  If we assume a given Schauder basis $\{ \bff_n \}_{n \in {\mathbb N}}$ in $Y$ or if we define a corresponding Schauder basis by setting $\bff_n = \bfq_n + \bfq_n^c$ then similar arguments can be used to find $Q \in {\mathcal B}(Y,Y_s)$ and $Q^c \in {\mathcal B}(Y,Y_r)$.  In practice the projection operators would be represented as infinite matrices relative to the respective Schauder bases.  The validity of the necessary elementary row and column operations required for solution of the various equations would depend on the properties of the Schauder bases.  It was shown by James \cite{jam1} that a Banach space with a Schauder basis is {\em reflexive} if and only if the basis is both {\em shrinking} and {\em boundedly complete}.  In such cases the elementary row and column operations would be justified.  These operations would also be justified if the Schauder basis is {\em unconditional}.  For an expanded discussion of bases in Banach space see Heil \cite[Chapter 3]{hei1}.

\section{Illustration of the main results}
\label{exmr}

We illustrate the main results with a pertinent example.

\begin{example}
\label{ex3}
{\rm  Let $p \in {\mathbb R}$ with $1 \leq p < \infty$ and let $X = Y = \ell^p$.  For each $z \in {\mathbb C}$ define $A(z) = A_0 + A_1z \in {\mathcal B}(\ell^p)$ by setting
$$
A_0 = \left[ \begin{array}{ccccccc}
1 & \alpha_1 & 0 & 0 & 0 & 0 & \cdots \\
0 & \beta & \alpha_2 & 0 & 0 &  0 & \cdots \\
0 & 0 & 1 & \alpha_3 & 0 & 0 & \cdots \\
0 & 0 & 0 & 0 & \alpha_4 & 0 & \cdots \\
0 & 0 & 0 & 0 & 1 & \alpha_5 & \cdots \\
0 & 0 & 0 & 0 & 0 & 0 & \cdots \\
\vdots & \vdots & \vdots & \vdots & \vdots & \vdots & \ddots \end{array} \right] \in {\mathcal B}(\ell^p) \quad \mbox{and} \quad A_1 = \left[ \begin{array}{ccccccc}
0 & 0 & 0 & 0 & 0 & 0 & \cdots \\
0 & 1 & 0 & 0 & 0 & 0 & \cdots \\
0 & 0 & 0 & 0 & 0 & 0 & \cdots \\
0 & 0 & 0 & 1 & 0 & 0 & \cdots \\
0 & 0 & 0 & 0 & 0 & 0 & \cdots \\
0 & 0 & 0 & 0 & 0 & 1 & \cdots \\
\vdots & \vdots & \vdots & \vdots & \vdots & \vdots & \ddots \end{array} \right] \in {\mathcal B}(\ell^p)
$$
where $|\alpha_j|^{1/j} \rightarrow 0$ as $j \rightarrow \infty$ and $\beta \neq 0$.

{\bf The Laurent series near zero.}  We wish to solve the fundamental equations and find a Laurent series for $R(z) = A(z)^{-1}$ on a region ${\mathcal U}_{\, 0,r}$ where $r > 0$.  If $\{ \bfx_{-n}\}_{n \in {\mathbb N}}$ is a singular Jordan chain with $A_0 \bfx_{-n} + A_1 \bfx_{-n-1} = \bfzero$ for all $n \in {\mathbb N}$ and if we write $\bfx_{-n} = [x_{-n,j}]$ then the general solution is given by
\begin{equation}
\label{gsol1}
\bfx_{-1} = \left[\hspace{-2mm} \begin{array}{c}
s_1 \\
- (1/\alpha_1)s_1 \\
s_3 \\
- (1/\alpha_3)s_3 \\
s_5 \\
- (1/\alpha_5)s_5 \\
\vdots \end{array} \hspace{-2mm} \right]\!\!, \ \bfx_{-2} = \left[ \hspace{-2mm} \begin{array}{c}
- \beta s_1 + \alpha_1 \alpha_2 s_3 \\
(\beta/\alpha_1) s_1 - \alpha_2 s_3 \\
\alpha_3 \alpha_4 s_5 \\
- \alpha_4 s_5 \\
\alpha_5 \alpha_6 s_7 \\
- \alpha_6 s_7 \\
\vdots \end{array} \hspace{-2mm} \right]\!\!, \ \bfx_{-3} = \left[ \hspace{-2mm} \begin{array}{c}
\beta^2 s_1 - \alpha_1 \alpha_2 \beta s_3 + \alpha_1 \alpha_2 \alpha_3 \alpha_4 s_5 \\
- (\beta^2/\alpha_1) s_1  + \alpha_2 \beta s_3 - \alpha_2 \alpha_3 \alpha_4  s_5 \\
\alpha_3 \alpha_4 \alpha_5 \alpha_6 s_7 \\
- \alpha_4 \alpha_5 \alpha_6 s_7 \\
\alpha_5 \alpha_6 \alpha_7 \alpha_8 s_9 \\
- \alpha_6 \alpha_7 \alpha_8 s_9 \\
\vdots \end{array} \hspace{-2mm} \right]\!\!, \cdots
\end{equation}
and so on {\em ad infinitum}, where the initial parameters $s_1,s_3, s_5, \ldots$ satisfy $|s_1|^p + |s_3|^p + |s_5|^p + \cdots < \infty$.  However it is clear by looking at the first component that if $s_1 \neq 0$ and $s_{2n+1} = 0$ for all $n \in {\mathbb N}$ then we will have $\|\bfx_{-n}\|^{1/n} \geq |s_1|^{1/n} |\beta|^{1-1/n} \rightarrow |\beta|$ as $n \rightarrow \infty$.  This shows that the general solution will not always generate an infinite-length Jordan chain with $\| \bfx_{-n}\|^{1/n} \rightarrow 0$ as $n \rightarrow \infty$.  Nevertheless, by judicious choice of the parameters, we can find an infinite-dimensional set of mutually linearly independent vectors, each of which generates a Jordan chain with only a finite number of nonzero elements in which case it is clear that we will have $\| \bfx_{-n}\|^{1/n} \rightarrow 0$ as $n \rightarrow \infty$.  Thus we can find a basis for our desired subspace $X_s = \{ \bfx_{-1} \mid A_0 \bfx_{-n} + A_1 \bfx_{-n-1} = \bfzero\ \mbox{for all}\ n \in {\mathbb N}\ \mbox{with}\ \| \bfx_{-n} \|^{1/n} \rightarrow 0\ \mbox{as}\ n \rightarrow \infty\}^a$. 

For instance we can define $X_s \subseteq \ell^p$ as the closed subspace generated by the vectors
\begin{eqnarray*}
\bfp_1 & = & - [\alpha_1 \alpha_2 \alpha_3/\beta]\bfe_1 + [\alpha_2 \alpha_3/\beta] \bfe_2 - \alpha_3 \bfe_3 + \bfe_4 \\
\bfp_2 & = & [\alpha_1 \cdots \alpha_5/\beta^2]\bfe_1 - [\alpha_2 \cdots \alpha_5/\beta^2] \bfe_2 - \alpha_5 \bfe_5 + \bfe_6 \\
\bfp_3 & = & - [\alpha_1 \cdots \alpha_7/\beta^3]\bfe_1 + [\alpha_2 \cdots \alpha_7/\beta^3] \bfe_2 - \alpha_7 \bfe_7 + \bfe_8 \\
\vdots & = & \vdots
\end{eqnarray*}
and so on.  In general $\bfp_m = (-1)^m[\alpha_1 \cdots \alpha_{2m+1}/\beta^m]\bfe_1 + (-1)^{m-1}[\alpha_2 \cdots \alpha_{2m+1}/\beta^m] \bfe_2 - \alpha_{2m+1} \bfe_{2m+1} + \bfe_{2m+2}$ for each $m \in {\mathbb N}$.   For each $m \in {\mathbb N}$ the base element $\bfp_m$ is found by choosing the vector $\bfs_m = (s_{1,m},s_{3,m},\ldots)$ of parameters in such a way that when $\bfx_{-1} = \bfp_m$ we have $\bfx_{-n} = \bfzero$ for $n > m$.  Thus the associated Jordan chain has only $m$ nonzero elements.   It follows that all finite linear combinations of these base elements also generate singular Jordan chains with only a finite number of nonzero elements.  Hence $\| \bfx_{-n}\|^{1/n} \rightarrow 0$ as $n \rightarrow \infty$ for all such chains.  

We wish to find infinite-length singular Jordan chains for $B(w) = B_0 + B_1w$ where $B_0 = A_1$ and $B_1 = A_0$ on some region $w \in {\mathcal U}_{1/r,\infty}$ where $r > 0$.  If $\{ \bfu_{-n}\}_{n \in {\mathbb N}}$ is a singular Jordan chain for $B(w)$ then we must have $B_0 \bfu_{-n} + B_1 \bfu_{-n-1} = \bfzero$.  If we write $\bfu_{-n} = [u_{-n,j}]$ then the general solution is given by
\begin{equation}
\label{gsol2}
\bfu_{-1} = \left[ \hspace{-2mm} \begin{array}{c}
t_1^{(1)} \\
t_2^{(1)} \\
t_3^{(1)} \\
t_4^{(1)} \\
t_5^{(1)} \\
t_6^{(1)} \\
t_7^{(1)} \\
t_8^{(1)} \\
\vdots \end{array} \hspace{-2mm} \right], \ \bfu_{-2} = \left[ \hspace{-2mm} \begin{array}{c}
t_1^{(2)} \\
- t_1^{(2)}/\alpha_1 \\
- ( t_2^{(1)} - \beta t_1^{(2)}/\alpha_1)/\alpha_2 \\
( t_2^{(1)} - \beta t_1^{(2)}/\alpha_1)/ (\alpha_2 \alpha_3) \\
- t_4^{(1)}/\alpha_4 \\
t_4^{(1)}/(\alpha_4 \alpha_5) \\
- t_6^{(1)}/ \alpha_6 \\
t_6^{(1)}/(\alpha_6 \alpha_7) \\
\vdots \end{array} \hspace{-2mm} \right], \ \bfu_{-3} = \left[ \hspace{-2mm} \begin{array}{c}
t_1^{(3)} \\
- t_1^{(3)}/\alpha_1 \\
( t_1^{(2)} + \beta t_1^{(3)})/ (\alpha_1 \alpha_2) \\
- ( t_1^{(2)} + \beta t_1^{(3)})/ (\alpha_1 \alpha_2 \alpha_3) \\
( t_2^{(1)} - \beta t_1^{(2)}/\alpha_1)/ (\alpha_2 \alpha_3 \alpha_4) \\
- ( t_2^{(1)} - \beta t_1^{(2)}/\alpha_1)/ (\alpha_2 \alpha_3 \alpha_4 \alpha_5) \\
- t_4^{(1)}/ (\alpha_4 \alpha_5 \alpha_6) \\
t_4^{(1)}/(\alpha_4 \alpha_5 \alpha_6 \alpha_7) \\
\vdots \end{array} \hspace{-2mm} \right], \cdots
\end{equation}
and so on where $\{t_{j}^{(1)}\}_{j \in {\mathbb N}}$ and $\{t_1^{(n+1)}\}_{n \in {\mathbb N}}$ are arbitrary parameters.  In this most general form it is clear that we may well have $\| \bfu_{-n} \|^{1/n} \rightarrow \infty$ as $n \rightarrow \infty$.  However if we choose $t_{2j}^{(1)} = 0$ for all $j \in {\mathbb N}+1$ and also choose $t_1^{(2)} = - t_2^{(1)}/\beta$ and $t_1^{(n+1)} = - t_1^{(n)}/\beta$ for each $n \in {\mathbb N}+1$ then we can see that
$$
\bfu_{-1} = t_1^{(1)} \bfe_1 + t_2^{(1)} \bfe_2 + t_3^{(1)} \bfe_3 + t_5^{(1)} \bfe_5 + t_7^{(1)} \bfe_7 + \cdots
$$
with
$$
\bfu_{-n} = (-1)^{n-1}(t_2^{(1)}/\beta^{n-1})\bfe_1 + (-1)^n (t_2^{(1)}/(\alpha_1 \beta^{n-1}) \bfe_2
$$
for all $n \geq 2$.  Now $\| \bfu_{-n} \| = \left(|t_2^{(1)}|/\beta^n \right) [1 + 1/\alpha_1^p]^{1/p} = C/\beta^n$ for $n \geq 2$ and so $\| \bfu_{-n} \|^{1/n} \rightarrow 1/ |\beta|$ as $n \rightarrow \infty$.  Hence $\{ \bfu_{-n} \}_{n \in {\mathbb N}}$ is an infinite-length singular Jordan chain for $B(w)$ on a region $w \in {\mathcal U}_{1/|\beta|, \infty} \Leftrightarrow z \in {\mathcal U}_{\, 0, |\beta|}$.  Thus $X_r$ can be defined as the closed subspace generated by the vectors
$$
\bfp_1^c = \bfe_1, \ \bfp_2^c = \bfe_2, \ \bfp_3^c = \bfe_3, \ \bfp_4^c = \bfe_5, \ \bfp_5^c = \bfe_7, \ \bfp_6^c = \bfe_9,\ \cdots
$$ 
and so on.  We found infinite-length singular Jordan chains $\{ \bfx_{-n}\}_{n \in {\mathbb N}}$ for $A(z)$ with $s = \lim_{n \rightarrow \infty} \| \bfx_{-n} \|^{1/n} = 0$ and $\{\bfu_{-n}\}_{n \in {\mathbb N}}$ for $B(w)$ with $1/r = \lim_{n \rightarrow \infty} \| \bfu_{-n} \|^{1/n} = 1/|\beta|$.  Since our analysis shows that there are no singular Jordan chains $\{ \bfu_{-n} \}_{n \in {\mathbb N}}$ for $B(w)$ with $\lim_{n \rightarrow \infty} \| \bfu_{-n} \|^{1/n} > 1/|\beta|$ it follows that  $A(z)$ is analytic for $z \in {\mathcal U}_{\, 0,|\beta|}$.   

To show that $\ell^p = X_s \oplus X_r$ it is sufficient to show that the equation
\begin{eqnarray*}
\lefteqn{x_1 \bfe_1 + x_2 \bfe_2 + x_3 \bfe_3 + x_4 \bfe_4 + \cdots + x_{2m} \bfe_{2m}} \\
& = & \xi_1 \bfp_1 + \xi_2 \bfp_2 + \cdots + \xi_{m-1} \bfp_{m-1} + \eta_1 \bfp_1^c + \eta_2 \bfp_2^c + \cdots + \eta_{m+1} \bfp_{m+1}^c 
\end{eqnarray*}
has a unique solution for each $m \in {\mathbb N}$.  By equating coefficients of $\bfe_{2k}$ and $\bfe_{2k-1}$ we can easily see that
$$
\xi_{k-1} = x_{2k}, \quad \eta_{k+1} = x_{2k-1} + \alpha_{2k-1} x_{2k}
$$
for each $k = 2,3,\ldots,m$.  Finally we obtain
\begin{eqnarray*}
\eta_1 & = & x_1 + [\alpha_1\alpha_2 \alpha_3/\beta] x_4 + (-1)[\alpha_1 \cdots \alpha_5/\beta^2] x_6 + \cdots + (-1)^{m-2}[\alpha_1 \cdots \alpha_{2m-1}/\beta^{m-1}] x_{2m}, \\
\eta_2 & = & x_2 + (-1)[\alpha_2 \alpha_3/\beta] x_4 + (-1)^2[\alpha_2 \cdots \alpha_5/\beta^2] x_6 + \cdots + (-1)^{m-1} [\alpha_2 \cdots \alpha_{2m-1}/\beta^{m-1}] x_{2m} 
\end{eqnarray*}
by equating coefficients of $\bfe_1$ and $\bfe_2$.  Thus the solution is unique.  Now it follows that
\begin{equation}
\label{psol}
P ( x_1 \bfe_1 + x_2 \bfe_2 + x_3 \bfe_3 + \cdots + x_{2m} \bfe_{2m}) =   x_4 \bfp_1 + x_6 \bfp_2 + \cdots + x_{2m} \bfp_{m-1}   
\end{equation}
and
\begin{eqnarray}
\label{pcsol}
\lefteqn{\hspace{-0.5cm}P^c ( x_1 \bfe_1 + x_2 \bfe_2 + x_3 \bfe_3 + \cdots + x_{2m} \bfe_{2m})} \nonumber \\
& = & \left\{ x_1  + [\alpha_1\alpha_2 \alpha_3/\beta] x_4 + (-1)[\alpha_1 \cdots \alpha_5/\beta^2] x_6 + \cdots + (-1)^{m-2}[\alpha_1\cdots \alpha_{2m-1} /\beta^{m-1}]x_{2m} \right\} \bfp_1^c \nonumber \\
& & \hspace{-1cm} + \left\{ x_2 + (-1)[\alpha_2 \alpha_3/\beta] x_4 + (-1)^2[\alpha_2 \cdots \alpha_5/\beta^2] x_6 + \cdots + (-1)^{m-1}[\alpha_2\cdots \alpha_{2m-1} /\beta^{m-1}]x_{2m} \right\} \bfp_2^c \nonumber \\
& & + \hspace{1mm} (x_3 + \alpha_3 x_4) \bfp_3^c + (x_5 + \alpha_5 x_6) \bfp_4^c + \cdots + (x_{2m-1} + \alpha_{2m-1} x_{2m}) \bfp_{m+1}^c
\end{eqnarray}
for each $\bfx \in \ell^p$.  We can represent each of the projections $P$ and $P^c$ in matrix form by defining column $j$ as the vector coefficient of $x_j$ in the respective solutions (\ref{psol}) and (\ref{pcsol}) written in terms of the standard basis.  Thus we have
$$
P = \left[ \begin{array}{ccccccc}
0 & 0 & 0 & -\alpha_1 \alpha_2 \alpha_3/\beta & 0 &  \alpha_1 \cdots \alpha_5/\beta^2 & \cdots \\
0 & 0 & 0 & \alpha_2 \alpha_3/\beta & 0 & - \alpha_2 \cdots \alpha_5/\beta^2 & \cdots \\
0 & 0 & 0 & -\alpha_3 & 0 & 0 & \cdots \\
0 & 0 & 0 & 1 & 0 & 0 & \cdots \\
0 & 0 & 0 & 0 & 0 & -\alpha_5 & \cdots \\
0 & 0 & 0 & 0 & 0 & 1 & \cdots \\
\vdots & \vdots & \vdots & \vdots & \vdots & \vdots & \ddots \end{array} \right] \in {\mathcal B}(\ell^p)
$$
and
$$
P^c = \left[ \begin{array}{ccccccc}
1 & 0 & 0 & \alpha_1 \alpha_2 \alpha_3/\beta & 0 & -\alpha_1 \cdots \alpha_5/\beta^2 & \cdots \\
0 & 1 & 0 & -\alpha_2 \alpha_3/\beta & 0 & \alpha_2 \cdots \alpha_5/\beta^2 & \cdots \\
0 & 0 & 1 & \alpha_3 & 0 & 0 & \cdots \\
0 & 0 & 0 & 0 & 0 & 0 & \cdots \\
0 & 0 & 0 & 0 & 1 & \alpha_5 & \cdots \\
0 & 0 & 0 & 0 & 0 & 0 & \cdots \\
\vdots & \vdots & \vdots & \vdots & \vdots & \vdots & \ddots \end{array} \right] \in {\mathcal B}(\ell^p).
$$
If $\{\bfp_n\}_{n \in {\mathbb N}}$ is a basis for $X_s$ and $\{ \bfp_n^c\}_{n \in {\mathbb N}}$ is a basis for $X_r$ and we define $\bfq_n = A_1\bfp_n$ and $\bfq_n^c = A_0 \bfp_n^c$ for each $n \in {\mathbb N}$ then $\{\bfq_n\}_{n \in {\mathbb N}}$ is a basis for $Y_s$ and $\{ \bfq_n^c\}_{n \in {\mathbb N}}$ is a basis for $Y_r$.  Therefore we can define $Y_s$ as the closed subspace generated by the vectors
$$ 
\bfq_1 = [\alpha_2 \alpha_3/\beta]\bfe_2 + \bfe_4, \ \bfq_2 = -[\alpha_2 \cdots \alpha_5/\beta^2]\bfe_2 + \bfe_6, \ \bfq_3 = [\alpha_2 \cdots \alpha_7/ \beta^3] \bfe_2 + \bfe_8, \ \ldots
$$
and so on.  In general $\bfq_{m-1} = (-1)^{m-2} [\alpha_2 \cdots \alpha_{2m-1}/\beta^{m-1}]\bfe_2 + \bfe_{2m}$ for each $m \in {\mathbb N}+1$.  The subspace $Y_r$ is the closed subspace generated by the vectors
$$
\bfq_1^c = \bfe_1, \ \bfq_2^c = \alpha_1 \bfe_1 +\beta \bfe_2, \ \bfq_3^c = \alpha_2 \bfe_2 + \bfe_3, \ \bfq_4^c = \alpha_4 \bfe_4 + \bfe_5, \ \ldots
$$
and so on.  In general $\bfq_{m+2}^c =  \alpha_{2m}\bfe_{2m} + \bfe_{2m+1}$ for all $m \in {\mathbb N}$.

To show that $\ell^p = Y_s \oplus Y_r$ it is sufficient to show that the equation
\begin{eqnarray*}
\lefteqn{ y_1 \bfe_1 + y_2 \bfe_2 + y_3 \bfe_3 + \cdots + y_{2m+1} \bfe_{2m+1} } \\
& = & \mu_1 \bfq_1 + \mu_2 \bfq_2 + \cdots + \mu_{m-1} \bfq_{m-1} + \nu_1 \bfq_1^c + \nu_2 \bfq_2^c + \cdots + \nu_{m+2} \bfq_{m+2}^c 
\end{eqnarray*}
has a unique solution for each $m \in {\mathbb N}$.  By equating coefficients of $\bfe_{2k}$ and $\bfe_{2k+1}$ we can easily see that $\nu_{k+2} =  y_{2k+1}$ for each $k=1,\ldots,m$ and $\mu_{k-1} = y_{2k} - \alpha_{2k}y_{2k+1}$ for each $k = 2,\ldots,m$.  If we equate coefficients of $\bfe_1$ and $\bfe_2$, we get
\begin{eqnarray*}
y_1 \hspace{-2mm} & = & \hspace{-2mm} \nu_1 + \alpha_1 \nu_2 \\
y_2 \hspace{-2mm} & = & \hspace{-2mm} [\alpha_2 \alpha_3/\beta]\mu_1 - [\alpha_2 \cdots \alpha_5/\beta^2] \mu_2 + \cdots + (-1)^{m-2}[\alpha_2 \cdots \alpha_{2m-1}/\beta^{m-1}] \mu_{m-1} + \beta \nu_2 + \alpha_2 \nu_3.
\end{eqnarray*}
But $\nu_3 = y_3$ and $\mu_{k-1} = y_{2k} - \alpha_{2k}y_{2k+1}$ for each $k = 2,\ldots,m$.  Thus, finally, we obtain 
\begin{eqnarray*}
\nu_2 \hspace{-2mm} & = & \hspace{-2mm} [1/\beta](y_2 - \alpha_2y_3) - [\alpha_2 \alpha_3/\beta^2](y_4 - \alpha_4 y_5) + \cdots + (-1)^{m-1}[\alpha_2 \cdots \alpha_{2m-1}/\beta^m](y_{2m} - \alpha_{2m} y_{2m+1}) \\
\nu_1 \hspace{-2mm} & = & \hspace{-2mm} y_1 - [\alpha_1/\beta](y_2 - \alpha_2y_3) + [\alpha_1 \alpha_2 \alpha_3/\beta^2](y_4 - \alpha_4 y_5) + \cdots + (-1)^m[\alpha_1 \cdots \alpha_{2m-1}/\beta^m](y_{2m} - \alpha_{2m} y_{2m+1}).
\end{eqnarray*}
Therefore the solution is unique.  Now it follows that
\begin{equation}
\label{qsol}
Q ( y_1 \bfe_1 + y_2 \bfe_2 + \cdots + y_{2m+1} \bfe_{2m+1}) =  (y_4 - \alpha_4 y_5)\bfq_1 + (y_6 - \alpha_6 y_7) \bfq_2 + \cdots + (y_{2m} - \alpha_{2m} y_{2m+1})  \bfq_{m-1}   
\end{equation}
and
\begin{eqnarray}
\label{qcsol}
\lefteqn{Q^c ( y_1 \bfe_1 + y_2 \bfe_2 + \cdots + y_{2m+1} \bfe_{2m+1})} \nonumber \\
& = & \left\{ y_1 - [\alpha_1/\beta](y_2 - \alpha_2y_3) + [\alpha_1 \alpha_2 \alpha_3/\beta^2](y_4 - \alpha_4 y_5) + \cdots \right. \nonumber \\
& & \hspace{7cm} \left. \cdots + (-1)^m[\alpha_1 \cdots \alpha_{2m-1}/\beta^m](y_{2m} - \alpha_{2m} y_{2m+1}) \right\} \bfq_1^c \nonumber \\
& & \hspace{0.5cm} + \left\{ [1/\beta](y_2 - \alpha_2y_3) - [\alpha_2 \alpha_3/\beta^2](y_4 - \alpha_4 y_5) + \cdots \right. \nonumber \\
& & \hspace{7cm} \left. \cdots + (-1)^{m-1}[\alpha_2 \cdots \alpha_{2m-1}/\beta^m](y_{2m} - \alpha_{2m} y_{2m+1}) \right\} \bfq_2^c \nonumber \\
& & \hspace{0.5cm} + \hspace{1mm} y_3 \bfq_3^c + y_5 \bfq_4^c + \cdots + y_{2m-1} \bfq_{m+1}^c
\end{eqnarray}
for each $\bfy \in \ell^p$.  We can represent each of the projections $Q$ and $Q^c$ in matrix form by defining column $j$ as the vector coefficient of $y_j$ in the respective solutions (\ref{qsol}) and (\ref{qcsol}) written in terms of the standard basis.  Thus we have
$$
Q = \left[ \begin{array}{cccccccc}
\hspace{0.25cm}0\hspace{0.25cm} & \hspace{0.25cm}0\hspace{0.25cm} & \hspace{0.25cm}0\hspace{0.25cm} & 0 & 0 & 0 & 0 & \cdots \\
0 & 0 & 0 & \alpha_2 \alpha_3/\beta & - \alpha_2 \alpha_3 \alpha_4/\beta & - \alpha_2 \cdots \alpha_5/\beta^2 &  \alpha_2 \cdots \alpha_6/\beta^2 &\cdots \\
0 & 0 & 0 & 0 & 0 & 0 & 0 & \cdots \\
0 & 0 & 0 & 1 & -\alpha_4 & 0 & 0 & \cdots \\
0 & 0 & 0 & 0 & 0 & 0 & 0 & \cdots \\
0 & 0 & 0 & 0 & 0 & 1 & -\alpha_6 & \cdots \\
0 & 0 & 0 & 0 & 0 & 0 & 0 & \cdots \\
\vdots & \vdots & \vdots & \vdots & \vdots & \vdots & \vdots & \ddots \end{array} \right] \in {\mathcal B}(\ell^p)
$$
and
$$
Q^c = \left[ \begin{array}{cccccccc}
\hspace{0.25cm}1\hspace{0.25cm} & \hspace{0.25cm}0\hspace{0.25cm} & \hspace{0.25cm}0\hspace{0.25cm} & 0 & 0 & 0 & 0 & \cdots \\
0 & 1 & 0 & - \alpha_2 \alpha_3/\beta & \alpha_2 \alpha_3\alpha_4/\beta & \alpha_2 \cdots \alpha_5/\beta^2 & - \alpha_2 \cdots \alpha_6/\beta^2 & \cdots \\
0 & 0 & 1 & 0 & 0 & 0 & 0 &\cdots \\
0 & 0 & 0 & 0 & \alpha_4 & 0 & 0 &\cdots \\
0 & 0 & 0 & 0 & 1 & 0 & 0 &\cdots \\
0 & 0 & 0 & 0 & 0 & 0 & \alpha_6 &\cdots \\
0 & 0 & 0 & 0 & 0 & 0 & 1 &\cdots \\
\vdots & \vdots & \vdots & \vdots & \vdots & \vdots & \vdots & \ddots \end{array} \right] \in {\mathcal B}(\ell^p).
$$
It was shown in \cite{alb2} that the basic solution $\{R_{-1}, R_0 \} \in {\mathcal B}(Y,X)$ is uniquely determined by the equations $R_{-1}A_1 = P$, $R_0A_0 = P^c$, $A_1R_{-1} = Q$, $A_0R_0 = Q^c$, $R_{-1} = PR_{-1}Q$ and $R_0 = P^cR_0Q^c$.  We have no general solution procedure to recommend for solving these equations once $P$ and $Q$ are known but in this particular instance we can proceed as follows.  Note that $\{R_{-1}, R_0 \} \in {\mathcal B}(\ell^p)$.  Solve the equations $R_{-1}A_1 = P$ and $R_0A_0 = P^c$ to find a partial solution for $R_{-1}$ and a complete solution for $R_0$.  Then solve the equations $A_1R_{-1} = Q$ and $R_{-1} = P R_{-1} Q$ to complete the determination of $R_{-1}$.  We obtain the solutions
$$
R_0 = \left[ \begin{array}{cccccccc}
1 & -\alpha_1/\beta & \alpha_1\alpha_2/\beta & \alpha_1 \alpha_2 \alpha_3/\beta^2 & - \alpha_1 \cdots \alpha_4/\beta^2 & -\alpha_1 \cdots \alpha_5/\beta^3 & \alpha_1 \cdots \alpha_6/\beta^3 & \cdots \\
0 & 1/\beta & - \alpha_2/\beta & - \alpha_2 \alpha_3/\beta^2 &  \alpha_2 \alpha_3 \alpha_4/\beta^2 & \alpha_2 \cdots \alpha_5/ \beta^3 &  - \alpha_2 \cdots \alpha_6/\beta^3 & \cdots \\
0 & 0 & 1 & 0 & 0 & 0 & 0 & \cdots \\
0 & 0 & 0 & 0 & 0 & 0 & 0 & \cdots \\
0 & 0 & 0 & 0 & 1 & 0 & 0 & \cdots \\
0 & 0 & 0 & 0 & 0 & 0 & 0 & \cdots \\
0 & 0 & 0 & 0 & 0 & 0 &  1 &\cdots \\
\vdots & \vdots & \vdots & \vdots & \vdots & \vdots & \vdots & \ddots \end{array} \right].
$$
and
$$
R_{-1} = \left[ \begin{array}{cccccccc}
\hspace{0.3cm}0\hspace{0.3cm} & \hspace{0.3cm}0\hspace{0.3cm} & \hspace{0.3cm}0\hspace{0.3cm} & -\alpha_1 \alpha_2 \alpha_3/\beta & \alpha_1 \cdots \alpha_4/\beta & \alpha_1 \cdots \alpha_5/\beta^2 & - \alpha_1 \cdots \alpha_6/\beta^2 & \cdots \\
0 & 0 & 0 & \alpha_2 \alpha_3/\beta & - \alpha_2 \alpha_3 \alpha_4/\beta & -\alpha_2 \cdots \alpha_5/\beta^2 &  \alpha_2 \cdots \alpha_6/\beta^2 & \cdots \\
0 & 0 & 0 & -\alpha_3 & \alpha_3  \alpha_4 & 0 & 0 & \cdots \\
0 & 0 & 0 & 1 & -\alpha_4 & 0 & 0 & \cdots \\
0 & 0 & 0 & 0 & 0 & -\alpha_5 & \alpha_5\alpha_6 & \cdots \\
0 & 0 & 0 & 0 & 0 & 1 & -\alpha_6 & \cdots \\
0 & 0 & 0 & 0 & 0 & 0 & 0 & \cdots \\
\vdots & \vdots & \vdots & \vdots & \vdots & \vdots & \vdots & \ddots \end{array} \right].
$$   
The entire sequence of Laurent series coefficients $\{ R_j\}_{j \in {\mathbb Z}}$ for the region ${\mathcal U}_{\, 0, |\beta|}$ can now be computed using the basic solution $\{R_{-1}, R_0 \}$ and the formul{\ae} $R_{-k} = (-1)^{k-1}(R_{-1}A_0)^{k-1}R_{-1}$ for each $k \in {\mathbb N}$ and $R_{\ell} = (-1)^{\ell} (R_0A_1)^{\ell} R_0$ for each $\ell \in {\mathbb N}-1$. 

{\bf The Laurent series near infinity.}  We wish to solve the fundamental equations and find a Laurent series for $R(z) = A(z)^{-1}$ on a region ${\mathcal U}_{\, s,\infty}$ where $s \geq |\beta|$.  We use a similar solution procedure to the one used to find the Laurent series on the region $0 < |z| < \beta$ but we omit some of the detailed calculations.  If $\{ \bfx_{-n}\}_{n \in {\mathbb N}}$ is a Jordan chain with $A_0 \bfx_{-n} + A_1 \bfx_{-n-1} = \bfzero$ for all $n \in {\mathbb N}$ and we set $s_{2k-1} = 1$ for some $k \in {\mathbb N}$ and assume that $s_{2j-1} = 0$ for all $j \in {\mathbb N}$ with $j \neq k$ in the general solution (\ref{gsol1}) obtained earlier then
$$
\bfx_{-1} = \bfe_{2k-1} - (1/\alpha_{2k-1}) \bfe_{2k}
$$
and for $n \geq 2k-1$ we have
$$
\bfx_{-n} = (-1)^{n-1}(\alpha_1 \alpha_2 \cdots \alpha_{2k-2}) \beta^{n-k} \bfe_1 + (-1)^n(\alpha_2 \cdots \alpha_{2k-2}) \beta^{n-k} \bfe_2.
$$ 
Since $\|\bfx_{-n}\| = |\alpha_2| \cdots |\alpha_{2k-2}| [\alpha_1^p  +1]^{1/p} \cdot |\beta|^{n-k}$ for $n \geq k$ it follows that $\| \bfx_{-n}\|^{1/n} \rightarrow |\beta|$ as $n \rightarrow \infty$.  If we define the subspace $X_s$ as the closure of the subspace spanned by the vectors
$$
\bfp_1 = \bfe_1 - (1/\alpha_1) \bfe_2, \quad  \bfp_2 = \bfe_3 - (1/\alpha_3) \bfe_4, \quad  \bfp_3 = \bfe_5 - (1/\alpha_5) \bfe_6, \quad \ldots
$$
then $X_s$ is the closed subspace defined by the points that generate all infinite-length singular Jordan chains for $A(z)$ on some region $z \in {\mathcal U}_{\, |\beta|,r}$ where $r > |\beta|$.  In similar fashion if $\{\bfu_{-n}\}_{n \in {\mathbb N}}$ is an infinite-length Jordan chain for $B(w)$ that satisfies the equations $B_0 \bfu_{-n} + B_1 \bfu_{-n-1} = \bfzero$ for all $n \in {\mathbb N}$ and we choose $t_{2k}^{(1)} = 0$ for all $k \in {\mathbb N}$ and $t_1^{(n+1)} = 0$ for all $n \in {\mathbb N}$ in the general solution (\ref{gsol2}) obtained earlier then we have $\bfu_{-n} = \bfzero$ for all $n \geq 2$ and hence $\| \bfu_{-n} \|^{1/n} \rightarrow 0$ as $n \rightarrow \infty$.  Thus the space $X_r$ defined by the closure of the linear space generated by the infinite set of vectors
$$
\bfp_1^c =  \bfe_1, \quad \bfp_2^c = \bfe_3, \quad \bfp_3^c = \bfe_5, \quad \ldots
$$
will generate all infinite-length singular Jordan chains for $B(w)$ on some region $w \in {\mathcal U}_{0, 1/s} \Leftrightarrow z \in {\mathcal U}_{s, \infty}$.

We found infinite-length singular Jordan chains $\{ \bfx_{-n}\}_{n \in {\mathbb N}}$ for $A(z)$ with $s = \lim_{n \rightarrow \infty} \| \bfx_{-n} \|^{1/n} = |\beta|$ and $\{\bfu_{-n}\}_{n \in {\mathbb N}}$ for $B(w)$ with $1/r = \lim_{n \rightarrow \infty} \| \bfu_{-n} \|^{1/n} = 0$.  Since there are no other possible singular Jordan chains $\{\bfx_{-n}\}_{n \in {\mathbb N}}$ with $\lim_{n \rightarrow \infty} \| \bfx_{-n} \|^{1/n} > |\beta|$ it follows that $A(z)$ is analytic for $z \in {\mathcal U}_{\,|\beta|, \infty}$.

We can now show that $\ell^p = X_s \oplus X_r$ and hence deduce that the two relevant projections can be represented in matrix form as
$$
P = [ \bfzero, -\alpha_1\bfe_1 + \bfe_2, \bfzero, -\alpha_3 \bfe_3 + \bfe_4, \bfzero, -\alpha_5 \bfe_5 + \bfe_6, \cdots ] \in {\mathcal B}(\ell^p)
$$
and
$$
P^c = [\bfe_1, \alpha_1 \bfe_1, \bfe_3, \alpha_3 \bfe_3, \bfe_5, \alpha_5 \bfe_5,\ldots ] \in {\mathcal B}(\ell^p).
$$
If $\{\bfp_n\}_{n \in {\mathbb N}}$ is a basis for $X_s$ and $\{ \bfp_n^c\}_{n \in {\mathbb N}}$ is a basis for $X_r$ and we define $\bfq_n = A_1\bfp_n$ and $\bfq_n^c = A_0 \bfp_n^c$ for each $n \in {\mathbb N}$ then $\{\bfq_n\}_{n \in {\mathbb N}}$ is a basis for $Y_s$ and $\{ \bfq_n^c\}_{n \in {\mathbb N}}$ is a basis for $Y_r$.  Therefore we can define $Y_s$ as the closed subspace generated by the vectors
$$
\bfq_1 = -(1/\alpha_1)\bfe_2,\ \bfq_2 = -(1/\alpha_3)\bfe_4,\ \bfq_3 = -(1/\alpha_5)\bfe_6,\ \ldots,\ \ldots
$$
and so on.  In general $\bfq_m = -(1/\alpha_{2m-1})\bfe_{2m}$ for each $m \in {\mathbb N}$.  The subspace $Y_r$ is the closed subspace generated by the vectors
$$
\bfq_1^c = \bfe_1,\ \bfq_2^c = \alpha_2\bfe_2 + \bfe_3,\ \bfq_3^c = \alpha_4\bfe_4 + \bfe_5,\ \ldots
$$
and so on.  In general $\bfq_{m+1}^c = \alpha_{2m}\bfe_{2m} + \bfe_{2m+1}$ for each $m \in {\mathbb N}$.

Now we can show that $\ell^p = Y_s \oplus Y_r$ and hence represent the two relevant projections in matrix form as
$$
Q = [\bfzero, \bfe_2, - \alpha_2 \bfe_2, \bfe_4, -\alpha_4 \bfe_4, \bfe_6, -\alpha_6 \bfe_6, \ldots ] \in {\mathcal B}(\ell^p)
$$
and
$$
Q^c = [ \bfe_1, \bfzero, \alpha_2\bfe_2 + \bfe_3, \bfzero, \alpha_4 \bfe_4 + \bfe_5, \bfzero, \alpha_6 \bfe_6 + \bfe_7, \ldots ] \in {\mathcal B}(\ell^p).
$$
To find the basic solution $\{R_{-1}, R_0 \} \in {\mathcal B}(\ell^p)$ we must solve the equations $R_{-1}A_1 = P$, $R_0A_0 = P^c$, $A_1R_{-1} = Q$, $A_0R_0 = Q^c$, $R_{-1} = PR_{-1}Q$ and $R_0 = P^cR_0Q^c$.   We proceed as follows.  We solve the equations $R_{-1}A_1 = P$ and $R_0A_0 = P^c$ to find a partial solution to $R_{-1}$ and a complete solution to $R_0$.  We then solve $A_1R_{-1} = Q$ and $R_{-1} = P R_{-1} Q$ to complete the solution process.  We obtain
$$
R_{-1} = \left[ \begin{array}{cccccccc}
0 & -\alpha_1 & \alpha_1\alpha_2 & 0 & 0 &  0 & 0 & \cdots \\
0 & 1 & -\alpha_2 & 0 & 0 & 0 & 0 & \cdots \\
0 & 0 & 0 & -\alpha_3 & \alpha_3 \alpha_4 & 0 & 0 & \cdots \\
0 & 0 & 0 & 1 & -\alpha_4 & 0 & 0 & \cdots \\
0 & 0 & 0 & 0 & 0 & -\alpha_5 & \alpha_5 \alpha_6 & \cdots \\
0 & 0 & 0 & 0 & 0 & 1 & -\alpha_6 & \cdots \\
\vdots & \vdots & \vdots & \vdots & \vdots & \vdots & \vdots & \ddots \end{array} \right] \ \mbox{and}\ \ R_0 = \left[ \begin{array}{ccccccc}
1 & 0 & 0 & 0 & 0 & 0 & \cdots \\
0 & 0 & 0 & 0 & 0 & 0 & \cdots \\
0 & 0 & 1 & 0 & 0 & 0 & \cdots \\
0 & 0 & 0 & 0 & 0 & 0 & \cdots \\
0 & 0 & 0 & 0 & 1 & 0 & \cdots \\
0 & 0 & 0 & 0 & 0 & 0 & \cdots \\
\vdots & \vdots & \vdots & \vdots & \vdots & \vdots & \ddots \end{array} \right].
$$
The entire sequence of Laurent series coefficients $\{ R_j\}_{j \in {\mathbb Z}}$ for the region ${\mathcal U}_{\, 0, |\beta|}$ can now be computed using the basic solution $\{R_{-1}, R_0 \}$ and the formul{\ae} $R_{-k} = (-1)^{k-1}(R_{-1}A_0)^{k-1}R_{-1}$ for each $k \in {\mathbb N}$ and $R_{\ell} = (-1)^{\ell} (R_0A_1)^{\ell} R_0$ for each $\ell \in {\mathbb N}-1$.  } $\hfill \Box$
\end{example}

\section{Conclusions and future research}
\label{cfr}

We suppose that $X$ and $Y$ are complex Banach spaces and that the generalized resolvent $R(z) = A(z)^{-1} \in {\mathcal B}(Y,X)$ of the linear pencil $A(z) = A_0 + A_1z \in {\mathcal B}(X,Y)$ has an essential singularity at the origin and is analytic on an annular region $z \in {\mathcal U}_{\, s,r}$ of the complex plane centred at the origin.  Hence the spectral set has a bounded component inside the annulus and an unbounded component outside the annulus.  In this case we found a complete characterization of the key spectral separation projection $P \in {\mathcal B}(X)$ in terms of the generating subspaces $X_s$ and $X_r$ for the corresponding infinite-length Jordan chains.  In particular we showed that in the domain space $P(X) = X_s$, $P^c(X) = X_r$ and $X = X_s \oplus X_r$.  We also established a corresponding structure in the range space for the projection $Q \in {\mathcal B}(Y)$ with $Q(Y) = Y_s$, $Q^c(X) = Y_r$ and $Y = Y_s \oplus Y_r$.  Our results are easily extended to cases where the essential singularity occurs at some point $z = a \neq 0$ and more generally to cases where the resolvent has a finite number of isolated essential singularities.  We showed in a particular case how the generating subspaces could be used to determine $P \in {\mathcal B}(X)$ and hence find a basic solution $\{R_{-1}, R_0 \} \in {\mathcal B}(Y,X)$ to the fundamental equations (\ref{lfe}) and (\ref{rfe}) on each annular region where the resolvent is analytic.  Our future research in this area will continue to seek improved systematic procedures for solving the fundamental equations in general Banach space.  

\section*{Acknowledgements}
\label{ack}

This research was funded by ARC Discovery Grant \#DP160101236  held by Phil Howlett, Amie Albrecht, Jerzy Filar and Konstantin Avrachenkov.  Geetika Verma is employed by the Grant as a Research Associate in the Scheduling and Control Group at the University of South Australia.


\begin{thebibliography}{00}

\bibitem{alb1} Amie R.~Albrecht, Phil G.~Howlett, Charles E.M.~Pearce, Necessary and sufficient conditions for the inversion of linearly-perturbed bounded linear operators on Banach space using Laurent series, J. Math. Anal. Appl. 383, (2011), 95\textendash 110.

\bibitem{alb2} Amie Albrecht, Phil Howlett, Charles Pearce, The fundamental equations for inversion of operator pencils on Banach space, J. Math. Anal. Appl. 413, (2014), 411\textendash 421. 

\bibitem{alb3} Amie Albrecht, Phil Howlett, Geetika Verma, Inversion of operator pencils on Hilbert space, J. Austral. Math. Soc., (2018), DOI: https://doi.org/10.1017/S1446788718000411.  Published online 21/12/2018.

\bibitem{alb4} Amie Albrecht, Phil Howlett, Geetika Verma, The fundamental equations for the generalized resolvent of an elementary pencil in a unital Banach algebra, Linear Algebra Appl., 574, (2019), 216\textendash 251.

\bibitem{avr2} K.E.~Avrachenkov, M.~Haviv, P.G.~Howlett, Inversion of analytic matrix functions that are singular at the origin, SIAM J. Matrix Anal. Appl. 22, (2001), 1175\textendash 1189.

\bibitem{avr3} Konstantin E. Avrachenkov, Jerzy A. Filar, Phil G. Howlett, Analytic Perturbation Theory and Its Applications OT 135, SIAM, Philadelphia, 2013.

\bibitem{bar1} Harm Bart, David C. Lay,  Poles of a generalised resolvent operator, P. Roy. Irish Acad. A. 74, (1974), 147\textendash 168. 

\bibitem{cas1} Pete Casazza, Ole Christensen, Diana T.~Stoeva, Frame expansions in separable Banach space, J. Math. Anal. Appl. 307, (2005), 710\textendash 723.

\bibitem{fra1} M. Franchi and P. Paruolo, Inverting a matrix function around a singularity via local rank factorization, SIAM J. Matix Anal. Appl. 37, (2016), 710\textendash 723. 

\bibitem{goh1} Israel Gohberg, Seymour Goldberg, Marinus A.~Kaashoek, Classes of Linear Operators Vol.~$1$, Operator Theory: Advances and Applications 49, Birkhauser Verlag, 1990.

\bibitem{hei1} Christopher Heil, A basis theory primer, School of Mathematics, Georgia Institute of Technology, Atlanta, Georgia, 1998, http://www.math.gatech.edu/$\sim$heil.

\bibitem{how1} P.G.~Howlett, Input retrieval in finite dimensional linear systems, ANZIAM J. (formerly J. Austral. Math. Soc. B) 23, (1982), 357\textendash 382.

\bibitem{how2} Phil Howlett, Konstantin Avrachenkov, Charles Pearce, Vladimir Ejov,  Inversion of analytically perturbed linear operators that are singular at the origin, J. Math. Anal. Appl. 353, (2009), 68\textendash 84.

\bibitem{how3} Phil Howlett, Amie Albrecht, Charles Pearce, Laurent series for inversion of linearly perturbed bounded linear operators on Banach space, J. Math. Anal. Appl. 366, (2010), 112\textendash 123.

\bibitem{jam1} Robert.C.~James, Bases and reflexivity of Banach spaces, Ann.~of Math. (2) 52,  (1950), 518\textendash 527.

\bibitem{kat1} Tosio Kato, Perturbation Theory for Linear Operators, Classics in Mathematics, Springer, 1995.

\bibitem{lan1} C.E.~Langenhop, The Laurent expansion for a nearly singular matrix, Linear Algebra Appl. 4, (1971), 329\textendash 340.
 
\bibitem{lan2} Carl E.~Langenhop,  On the Invertibility of a Nearly Singular Matrix, Linear Algebra Appl. 7, (1973), 361\textendash 365.

\bibitem{ros1}  N.J.~Rose, The Laurent expansion of a generalized resolvent with some applications, SIAM J. Math. Anal. 9, (1978), 751\textendash 758.

\bibitem{rot1} U. G. Rothblum, Resolvent expansion of matrices and applications, Linear Algebra Appl. 38, (1981), 33\textendash 49.

\bibitem{sch1} P.~Schweitzer, G.W.~Stewart, The Laurent expansion of pencils that are singular at the origin, Linear Algebra Appl. 183, (1993), 237\textendash 254.

\bibitem{stu1} F.~Stummel,  Diskrete Konvergenz Linearer Operatoren II, Math. Z. 120, (1971), 231\textendash 264 (in German).

\bibitem{wil1} Jon Wilkening,  An algorithm for computing Jordan chains and inverting analytic matrix functions, Linear Algebra Appl. 427, (2007), 6\textendash 25.

\bibitem{yos1} K\^{o}saku Yosida, Functional Analysis, Fifth Edition, Classics in Mathematics, Springer-Verlag, 1978.

\end{thebibliography}
\end{document}